\documentclass[a4paper,12pt]{article}
\usepackage[cp1251]{inputenc} % Для WiN-кодировки в MiKTeX
\usepackage[english]{babel}
\usepackage{amsfonts,amssymb}
\title{The homotopy invariance of cyclic homology of $\A$-algebras over rings.}
\author{S.V. Lapin}
\date{}
\newcommand{\F}{F_{\infty}}
\newcommand{\D}{D_{\infty}}
\newcommand{\bu}{\bullet}
\newcommand{\E}{E_{\infty}}
\newcommand{\A}{A_{\infty}}

\newcommand{\p}{\partial}

\begin{document}
\maketitle

\begin{abstract} In the present paper the cyclic homology functor
from the category of $A_\infty$-al\-geb\-ras over any commutative
unital ring $K$ to the category of graded $K$-modules is
constructed. Further, it is showed that this functor sends
homotopy equivalences of $A_\infty$-al\-geb\-ras into isomorphisms
of graded modules. As a corollary, it is obtained that the cyclic
homology of an $A_\infty$-algebra over any field is isomorphic to
the cyclic homology of the $A_\infty$-algebra of homologies for
the source $A_\infty$-algebra.
\end{abstract}

In \cite{Lapin}, on the basis of the combinatorial and homotopy
technique of differential modules with $\infty$-simplicial faces
\cite{Lap1}-\cite{Lap7} and $\D$-differential modules
\cite{Lap9}-\cite{Lap17} the cyclic bicomplex of an $\A$-algebra
over any commutative unital ring was constructed. This bicomplex
generalizes the cyclic bicomplex \cite{LQ} of an associative
algebra given over an arbitrary commutative unital ring. Further,
in \cite{Lapin}, cyclic homology of any $\A$-al\-geb\-ra over an
arbitrary commutative unital ring was defined as the homology of
the chain complex associated with the cyclic bicomplex of this
$\A$-algebra. The cyclic homology of $\A$-algebras over
commutative unital rings introduced in \cite{Lapin} generalizes
the cyclic homology of associative algebras over commutative
unital rings defined in \cite{LQ}. It is well known \cite{LQ} that
over fields of characteristic zero the cyclic homology introduced
in \cite{LQ} is isomorphic to the cyclic homology defined in
\cite{Con1} by using the complex of coinvariants for the action of
cyclic groups. Similar to this, in \cite{Lapin}, it was shown that
over fields of characteristic zero the cyclic homology of
$\A$-algebras introduced in \cite{Lapin} is isomorphic to the
cyclic homology of $\A$-algebras defined in \cite{PS} by using the
complex of coinvariants for the action of cyclic groups. Moreover,
in \cite{Lapin}, for homotopy unital $\A$-algebras over any
commutative unital rings, the analogue of the Connes–Tsygan exact
sequence was constructed.

On the other hand, in \cite{Kad}, it was shown that the structure
of an $\A$-algebra is homotopy invariant, i.e., the specified
structure is invariant under homotopy equivalences of differential
modules. In addition, in \cite{Kad}, it was established that the
homology of the $B$-con\-struc\-tion of an $\A$-al\-geb\-ra is
homotopy invariant, i.e., this homology is invariant under
homotopy equivalences of $\A$-al\-geb\-ras. Now note that the
cyclic bicomplex of an $\A$-al\-geb\-ra constructed in
\cite{Lapin} is the cyclic analogue of the $B$-con\-struc\-tion of
an $\A$-al\-geb\-ra, and the cyclic homology of an
$\A$-al\-geb\-ra is defined in \cite{Lapin} as the homology of
this cyclic analogue of the $B$-con\-struc\-tion. This gives rise
to an interesting natural question: do the cyclic homology of
$\A$-al\-geb\-ras is homotopy invariant under the homotopy
equivalences of $\A$-al\-geb\-ras? In present paper a positive
answer to this question is given.

The paper consists of three paragraphs. In the first paragraph, we
first recall necessary definitions related to the notion of a
cyclic module with $\infty$-simplicial faces or, more briefly, an
$C\F$-module \cite{Lapin}, which homotopy generalizes the notion
of a cyclic module with simplicial faces \cite{Con}. After that,
the category of $C\F$-modules is defined, namely, the notion of a
morphism of $C\F$-modules is introduced, and it is shown that the
composition of morphisms of $C\F$-modules is a morphism of
$C\F$-modules. Next, the concept of a homotopy between morphisms
of $C\F$-modules and the notion of a homotopy equivalence of
$C\F$-modules are introduced.

In the second paragraph, we first recall necessary definitions
related to the notion of a cyclic homology of $C\F$-modules
\cite{Lapin}. Next, it is shown that the cyclic homology of
$C\F$-modules defines the functor from the category of
$C\F$-modules to the category of graded modules. In addition, it
is shown that this functor sends homotopy equivalences of
$C\F$-modules into isomorphisms of graded modules.

In the third paragraph, we first recall necessary definitions
related to the notion of an $\A$-algebra \cite{Kad}. Next, we
recall the concept of a cyclic homology of $\A$-algebras over an
arbitrary commutative unital rings \cite{Lapin}. Then, by using
results of the second paragraph, it is shown that the cyclic
homology of $\A$-algebras defines the functor from the category of
$\A$-algebras to the category of graded modules. Moreover, it is
shown that this functor sends homotopy equivalences of
$\A$-algebras into isomorphisms of graded modules. As a corollary,
we obtain that the cyclic homology of an $\A$-algebra over any
field is isomorphic to the cyclic homology of the $\A$-algebra of
homologies for the source $\A$-algebra. In particular, it is
obtained that the cyclic homology of an associative differential
algebra over any field is isomorphic to the cyclic homology of the
$\A$-algebra of homologies for the source associative differential
algebra.

We proceed to precise definitions and statements. All modules and
maps of modules considered in this paper are, respectively,
$K$-modules and $K$-linear maps of modules, where $K$ is any
unital (i.e., with unit) commutative ring. \vspace{0.5cm}

\centerline{\bf \S\,1. Cyclic modules with $\infty$-simplicial
faces and} \centerline{\bf their morphisms and homotopies}
\vspace{0.5cm}

In what follows, by a bigraded module we mean any bigraded module
$X=\{X_{n,\,m}\}$, $n\geqslant 0$, $m\geqslant 0$, and by a
differential bigraded module, or, briefly, a differential module
$(X,d)$, we mean any bigraded module $X$ endowed with a
differential $d:X_{*,\bu}\to X_{*,\bu-1}$ of bidegree $(0,-1)$.

Recall that a differential module with simplicial faces is defined
as a differential module $(X,d)$ together with a family of module
maps $\p_i:X_{n,\bu}\to X_{n-1,\bu}$, $0\leqslant i\leqslant n$,
which are maps of differential modules and satisfy the simplicial
commutation relations $\p_i\p_j=\p_{j-1}\p_i$, $i<j$. The maps
$\p_i:X_{n,\bu}\to X_{n-1,\bu}$ are called the simplicial face
operators or, more briefly, the simplicial faces of the
differential module $(X,d)$.

Now, we recall the notion of a differential module with
$\infty$-simplicial faces \cite{Lap1} (see also
\cite{Lap2}-\cite{Lap7}), which is a homotopy invariant analogue
of the notion of a differential module with simplicial faces.

Let $\Sigma_k$ be the symmetric group of permutations on a
$k$-element set. Given an arbitrary permutation
$\sigma\in\Sigma_k$ and any $k$-tuple of nonnegative integers
$(i_1,\dots,i_k)$, where $i_1<\dots<i_k$, we consider the
$k$-tuple $(\sigma(i_1),\dots,\sigma(i_k))$, where $\sigma$ acts
on the $k$-tuple $(i_1,\dots,i_k)$ in the standard way, i.e.,
permutes its components. For the $k$-tuple
$(\sigma(i_1),\dots,\sigma(i_k))$, we define a $k$-tuple
$(\widehat{\sigma(i_1)},\dots,\widehat{\sigma(i_k)})$ by the
following formulae $$\widehat{\sigma
(i_s)}=\sigma(i_s)-\alpha(\sigma(i_s)),\quad 1\leqslant s\leqslant
k,$$ where each $\alpha(\sigma(i_s))$ is the number of those
elements of $(\sigma(i_1),\dots,\sigma(i_s),\dots\sigma(i_k))$ on
the right of $\sigma(i_s)$ that are smaller than $\sigma(i_s)$.

A differential module with $\infty$-simplicial faces or, more
briefly, an $\F$-module $(X,d,\p)$ is defined as a differential
module $(X,d)$ together with a family of module maps
$\p=\{\p_{(i_1,\dots ,i_k)}:X_{n,\bu}\to X_{n-k,\bu+k-1}\}$,
$1\leqslant k\leqslant n$, $0\leqslant i_1<\dots<i_k\leqslant n$,
$i_1,\dots,i_k\in\mathbb{Z}$, which satisfy the relations
$$d(\p_{(i_1,\dots,i_k)})=\sum_{\sigma\in\Sigma_k}\sum_{I_{\sigma}}
(-1)^{{\rm sign}(\sigma)+1}
\p_{(\widehat{\sigma(i_1)},\dots,\widehat{\sigma(i_m)})}\,
\p_{(\widehat{\sigma(i_{m+1})},\dots,\widehat{\sigma(i_k)})},\eqno(1.1)$$
where $I_\sigma$ is the set of all partitions of the $k$-tuple
$(\widehat{\sigma(i_1)},\dots,\widehat{\sigma(i_k)})$ into two
tuples $(\widehat{\sigma(i_1)},\dots,\widehat{\sigma(i_m)})$ and
$(\widehat{\sigma(i_{m+1})},\dots,\widehat{\sigma(i_k)})$,
$1\leqslant m\leqslant k-1$, such that the conditions
$\widehat{\sigma(i_1)}<\dots<\widehat{\sigma(i_m)}$ and
$\widehat{\sigma(i_{m+1})}<\dots<\widehat{\sigma(i_k)}$ holds.

The family of maps $\p=\{\p_{(i_1,\dots ,i_k)}\}$ is called the
$\F$-differential of the $\F$-module $(X,d,\widetilde{\p})$. The
maps $\p_{(i_1,\dots ,i_k)}$ that form the $\F$-differential of an
$\F$-module $(X,d,\p)$ are called the $\infty$-simplicial faces of
this $\F$-module.

It is easy to show that, for $k=1,2,3$, relations $(1.1)$ take,
respectively, the following view $$d(\p_{(i)})=0,\quad i\geqslant
0,\quad d(\p_{(i,j)})=\p_{(j-1)}\p_{(i)}-\p_{(i)}\p_{(j)},\quad
i<j,$$
$$d(\p_{(i_1,i_2,i_3)})=-\p_{(i_1)}\p_{(i_2,i_3)}-\p_{(i_1,i_2)}\p_{(i_3)}-
\p_{(i_3-2)}\p_{(i_1,i_2)}-$$
$$-\,\p_{(i_2-1,i_3-1)}\p_{(i_1)}+\p_{(i_2-1)}\p_{(i_1,i_3)}+\p_{(i_1,i_3-1)}\p_{(i_2)},
\quad i_1<i_2<i_3.$$

It is easy to check that, for any permutation $\sigma\in\Sigma_k$
and any $k$-tuple $(i_1,\dots,i_k)$, where $i_1<\dots<i_k$, the
conditions $\widehat{\sigma(i_1)}<\dots<\widehat{\sigma(i_m)}$ and
$\widehat{\sigma(i_{m+1})}<\dots<\widehat{\sigma(i_k)}$ are
equivalent to the conditions $\sigma(i_1)<\dots<\sigma(i_m)$ and
$\sigma(i_{m+1})<\dots<\sigma(i_k)$. This readily implies that the
$k$-tuple
$(\widehat{\sigma(i_{m+1})},\dots,\widehat{\sigma(i_k)})$, which
specified in $(1.1)$, coincides with the $k$-tuple
$(\sigma(i_{m+1}),\dots,\sigma(i_k))$.

Simplest examples of differential modules with $\infty$-simplicial
faces are differential modules with simplicial faces. Indeed,
given any differential module with simplicial faces $(X,d,\p_i)$,
we can define the $\F$-differential $\p=\{\p_{(i_1,\dots
,i_k)}\}:X\to X$ by setting $\p_{(i)}=\p_i$, $i\geqslant 0$, and
$\p_{(i_1,\dots,i_k)}=0$, $k>1$, thus obtaining the differential
module with $\infty$-simplicial faces $(X,d,\p)$.

It is worth mentioning that the notion of an differential module
with $\infty$-simplicial faces specified above is a part of the
general notion of a differential $\infty$-simplicial module
introduced in \cite{Lap3} by using the homotopy technique of
differential Lie modules over curved colored coalgebras.

Recall \cite{Con} that a cyclic differential module with
simplicial faces $(X,d,\p_i,t)$ is defined as a differential
module with simplicial faces $(X,d,\p_i)$ equipped with a family
of module maps $t=\{t_n:X_{n,\bu}\to X_{n,\bu}\}$, $n\geqslant 0$,
which satisfy the following relations: $$dt_n=t_nd,\quad
t_n^{\,n+1}=1_{X_{n,\bu}},\quad n\geqslant 0,$$
$$\p_it_n=t_{n-1}\p_{i-1},\quad 0<i\leqslant n,\quad
\p_0t_n=\p_n.$$

Now, let us recall \cite{Lapin} that a cyclic differential module
with $\infty$-simplicial faces or, more briefly, an $C\F$-module
$(X,d,\p,t)$ is defined as any $\F$-module $(X,d,\p)$ together
with a family of module maps $t=\{t_n:X_{n,\bu}\to X_{n,\bu}\}$,
$n\geqslant 0$, which satisfy the following relations:
$$dt_n=t_nd,\quad t_n^{\,n+1}=1_{X_{n,\bu}},\quad n\geqslant 0,$$
$$\p_{(i_1,\dots,i_k)}t_n=\left\{\begin{array}{ll}
t_{n-k}\p_{(i_1-1,\dots,i_k-1)},& i_1>0,\\
(-1)^{k-1}\p_{(i_2-1,\dots,i_k-1,n)},& i_1=0.\\
\end{array}\right.\eqno(1.2)$$

The family of maps $\p=\{\p_{(i_1,\dots ,i_k)}\}$ is called the
$\F$-differential of the $C\F$-mo\-du\-le $(X,d,\p,t)$. The maps
$\p_{(i_1,\dots ,i_k)}$ are called the $\infty$-simplicial faces
of this $C\F$-module.

Simplest examples of $C\F$-modules are cyclic differential modules
with simplicial faces. Indeed, given any cyclic differential
module with simplicial faces $(X,d,\p_i,t)$, we can define the
$\F$-differential $\p=\{\p_{(i_1,\dots ,i_k)}\}:X\to X$ by setting
$\p_{(i)}=\p_i$, $i\geqslant 0$, and $\p_{(i_1,\dots,i_k)}=0$,
$k>1$, thus obtaining the $C\F$-module $(X,d,\p,t)$.

It is worth mentioning that the notion of a $C\F$-module specified
above is a part of the general notion of a cyclic
$\infty$-simplicial module introduced in \cite{Lapin1} by using
the homotopy technique of differential modules over curved colored
coalgebras.

Now, we recall that a map $f:(X,d,\p_i)\to (Y,d,\p_i)$ of
differential modules with simplicial faces is defined as a map of
differential modules $f:(X,d)\to(Y,d)$ that satisfies the
relations $\p_if=f\p_i$, $i\geqslant 0$.

Let us consider the notion of a morphism of differential modules
with $\infty$-simplicial faces \cite{Lap1} (see also \cite{Lap4}),
which homotopically generalizes the notion of a map differential
modules with simplicial faces.

A morphism of $\F$-modules $f:(X,d,\p)\to (Y,d,\p)$ is defined as
a family of module maps $f=\{f_{(i_1,\dots,i_k)}:X_{n,\bu}\to
Y_{n-k,\bu+k}\}$, $0\leqslant k\leqslant n$, $0\leqslant
i_1<\dots<i_k\leqslant n$, $i_1,\dots,i_k\in\mathbb{Z}$, (at $k=0$
we will use the denotation $f_{(\,\,)}$), which satisfy the
relations
$$d(f_{(i_1,\dots,i_k)})=-\p_{(i_1,\dots,i_k)}f_{(\,\,)}+f_{(\,\,)}\p_{(i_1,\dots,i_k)}\,+$$
$$+\sum_{\sigma\in\Sigma_k}\sum_{I_{\sigma}}(-1)^{{\rm
sign}(\sigma)+1}\p_{(\widehat{\sigma(i_1)},\dots,
\widehat{\sigma(i_m)})}f_{(\widehat{\sigma(i_{m+1})},\dots,
\widehat{\sigma(i_k)})}\,-$$ $$-\,f_{(\widehat{\sigma(i_1)},\dots,
\widehat{\sigma(i_m)})}\p_{(\widehat{\sigma(i_{m+1})},\dots,
\widehat{\sigma(i_k)})},\eqno(1.3)$$ where $I_{\sigma}$ is the
same as in $(1.1)$. The maps $f_{(i_1,\dots ,i_k)}\in f$ are
called the components of the morphism $f:(X,d,\p)\to (Y,d,\p)$.

For example, at $k=0,1,2,3$ the relations $(1.3)$ take,
respectively, the following view $$d(f_{(\,\,)})=0,\qquad
d(f_{(i)})=f_{(\,\,)}\p_{(i)}-\p_{(i)}f_{(\,\,)},\quad i\geqslant
0,$$
$$d(f_{(i,j)})=-\p_{(i,j)}f_{(\,\,)}+f_{(\,\,)}\p_{(i,j)}-\p_{(i)}f_{(j)}+
\p_{(j-1)}f_{(i)}+ f_{(i)}\p_{(j)}-f_{(j-1)}\p_{(i)},\quad i<j,$$
$$d(f_{(i_1,i_2,i_3)})=-\p_{(i_1,i_2,i_3)}f_{(\,\,)}+f_{(\,\,)}\p_{(i_1,i_2,i_3)}-
\p_{(i_1)}f_{(i_2,i_3)}-\p_{(i_1,i_2)}f_{(i_3)}-
\p_{(i_3-2)}f_{(i_1,i_2)}-$$
$$-\,\p_{(i_2-1,i_3-1)}f_{(i_1)}+\p_{(i_2-1)}f_{(i_1,i_3)}+\p_{(i_1,i_3-1)}f_{(i_2)}+
f_{(i_1)}\p_{(i_2,i_3)}+f_{(i_1,i_2)}\p_{(i_3)}+$$
$$+f_{(i_3-2)}\p_{(i_1,i_2)}+
f_{(i_2-1,i_3-1)}\p_{(i_1)}-f_{(i_2-1)}\p_{(i_1,i_3)}-f_{(i_1,i_3-1)}\p_{(i_2)},\quad
i_1<i_2<i_3.$$

Now, we recall \cite{Lap1} that a composition of an arbitrary
given morphisms of $\F$-mo\-du\-les $f:(X,d,\p)\to (Y,d,\p)$ and
$g:(Y,d,\p)\to (Z,d,\p)$ is defined as  a morphism of $\F$-modules
$gf:(X,d,\p)\to (Z,d,\p)$ whose components are defined by
$$(gf)_{(i_1,\dots,i_k)}=\sum_{\sigma\in\Sigma_k}\sum_{I'_{\sigma}}
(-1)^{{\rm sign}(\sigma)}
g_{(\widehat{\sigma(i_1)},\dots,\widehat{\sigma(i_m)})}
f_{(\widehat{\sigma(i_{m+1})},\dots,\widehat{\sigma(i_k)})},\eqno(1.4)$$
where $I'_\sigma$ is the set of all partitions of the $k$-tuple
$(\widehat{\sigma(i_1)},\dots,\widehat{\sigma(i_k)})$ into two
tuples $(\widehat{\sigma(i_1)},\dots,\widehat{\sigma(i_m)})$ and
$(\widehat{\sigma(i_{m+1})},\dots,\widehat{\sigma(i_k)})$,
$0\leqslant m\leqslant k$, such that the conditions
$\widehat{\sigma(i_1)}<\dots<\widehat{\sigma(i_m)}$ and
$\widehat{\sigma(i_{m+1})}<\dots<\widehat{\sigma(i_k)}$ holds.

For example, at $k=0,1,2,3$ the formulae $(1.4)$ take,
respectively, the following form
$$(gf)_{(\,\,)}=g_{(\,\,)}f_{(\,\,)},\qquad(gf)_{(i)}=g_{(\,\,)}f_{(i)}+g_{(i)}f_{(\,\,)},$$
$$(gf)_{(i_1,i_2)}=g_{(\,\,)}f_{(i_1,i_2)}+g_{(i_1,i_2)}f_{(\,\,)}+
g_{(i_1)}f_{(i_2)}-g_{(i_2-1)}f_{(i_1)},\qquad i_1<i_2,$$
$$(gf)_{(i_1,i_2,i_3)}=g_{(\,\,)}f_{(i_1,i_2,i_3)}+g_{(i_1,i_2,i_3)}f_{(\,\,)}+
g_{(i_1)}f_{(i_2,i_3)}+g_{(i_1,i_2)}f_{(i_3)}+$$
$$+\,g_{(i_3-2)}f_{(i_1,i_2)}+
g_{(i_2-1,i_3-1)}f_{(i_1)}-g_{(i_2-1)}f_{(i_1,i_3)}-g_{(i_1,i_3-1)}f_{(i_2)},\quad
i_1<i_2<i_3.$$

{\bf Definition 1.1}. A morphism of $C\F$-modules $f:(X,d,\p,t)\to
(Y,d,\p,t)$ is defined as any morphism of $\F$-modules
$f:(X,d,\p)\to (Y,d,\p)$ whose components satisfy the following
conditions: $$f_{(\,\,)}t_n=t_nf_{(\,\,)},\quad
f_{(i_1,\dots,i_k)}t_n=\left\{\begin{array}{lll}
t_{n-k}f_{(i_1-1,\dots,i_k-1)},&k\geqslant 1,& i_1>0,\\
(-1)^{k-1}f_{(i_2-1,\dots,i_k-1,n)},&k\geqslant 1,& i_1=0.\\
\end{array}\right.\eqno(1.5)$$

By using the fact that any morphism of $C\F$-modules is a morphism
of $\F$-mo\-du\-les we define the composition of morphisms of
$C\F$-modules as a composition of morphisms of $\F$-modules.

{\bf Theorem 1.1}. The composition of morphisms of $C\F$-modules
is a morphism of $C\F$-modules.

{\bf Proof}. For an arbitrary morphisms of $C\F$-modules
$f:(X,d,\p,t)\to (Y,d,\p,t)$ and $g:(Y,d,\p,t)\to (Z,d,\p,t)$, we
need to check that components of the morphism of $\F$-modules
$gf:(X,d,\p)\to (Y,d,\p)$ satisfy the relations $(1.5)$. It is
clearly that at $k=0$ we have $(gf)_{(\,\,)}t_n=t_n(gf)_{(\,\,)}$.
Now note, for any $k$-tuple $(i_1,\dots,i_k)$ and any permutation
$\sigma\in\Sigma_k$, where $k\geqslant 1$ and $0<i_1<\dots<i_k$,
the $k$-tuple
$(\widehat{\sigma(i_1)},\dots,\widehat{\sigma(i_k)})$ satisfies
the conditions
$\widehat{\sigma(i_1)}>0,\dots,\widehat{\sigma(i_k)}>0$. By using
these conditions we obtain
$$(gf)_{(i_1,\dots,i_k)}t_n=\sum_{\sigma\in\Sigma_k}\sum_{I_{\sigma}}
(-1)^{{\rm sign}(\sigma)}
g_{(\widehat{\sigma(i_1)},\dots,\widehat{\sigma(i_m)})}
f_{(\widehat{\sigma(i_{m+1})},\dots,\widehat{\sigma(i_k)})}t_n=$$
$$=t_{n-k}\sum_{\sigma\in\Sigma_k}\sum_{I_{\sigma}} (-1)^{{\rm
sign}(\sigma)}
g_{(\widehat{\sigma(i_1)}-1,\dots,\widehat{\sigma(i_m)}-1)}
f_{(\widehat{\sigma(i_{m+1})}-1,\dots,\widehat{\sigma(i_k)}-1)}=$$
$$=t_{n-k}\sum_{\sigma\in\Sigma_k}\sum_{I_{\sigma}} (-1)^{{\rm
sign}(\sigma)}
g_{(\widehat{\sigma(i_1-1)},\dots,\widehat{\sigma(i_m-1)})}
f_{(\widehat{\sigma(i_{m+1}-1)},\dots,\widehat{\sigma(i_k)-1})}=t_{n-k}(gf)_{(i_1-1,\dots,i_k-1)}.$$
Now we show that at $k\geqslant 1$ and $i_1=0$ the relations
$(1.5)$ holds, namely, we show that the equality
$(gf)_{(0,i_2,\dots,i_k)}t_n=(-1)^{k-1}(gf)_{(i_2-1,\dots,i_k-1,n)}$
is true. By the definition of a composition we have the equalities
$$(gf)_{(0,i_2,\dots,i_k)}t_n=g_{(\,\,)}f_{(0,i_2,\dots,i_k)}t_n+g_{(0,i_2,\dots,i_k)}f_{(\,\,)}t_n+$$
$$+\sum_{\sigma\in\Sigma_k}\sum_{I'_{\sigma}} (-1)^{{\rm
sign}(\sigma)}
g_{(\widehat{\sigma(0)},\widehat{\sigma(i_2)},\dots,\widehat{\sigma(i_m)})}
f_{(\widehat{\sigma(i_{m+1})},\dots,\widehat{\sigma(i_k)})}t_n,\eqno(1.6)$$
$$(-1)^{k-1}(gf)_{(i_2-1,\dots,i_k-1,n)}=(-1)^{k-1}g_{(\,\,)}f_{(i_2-1,\dots,i_k-1,n)}+
(-1)^{k-1}g_{(i_2-1,\dots,i_k-1,n)}f_{(\,\,)}+$$
$$+(-1)^{k-1}\sum_{\varrho\in\Sigma_k}\sum_{I'_{\varrho}}
(-1)^{{\rm sign}(\varrho)}
g_{(\widehat{\varrho(i_2-1)},\dots,\widehat{\varrho(i_{m+1}-1)})}
f_{(\widehat{\varrho(i_{m+2}-1)},\dots,\widehat{\varrho(i_k-1)},\widehat{\varrho(n)})}.\eqno(1.7)$$
Let us show that each summand on the right-hand side of $(1.6)$ is
equal to some summand on the right-hand side of $(1.7)$. It is
easy to see that
$g_{(\,\,)}f_{(0,i_2,\dots,i_k)}t_n=(-1)^{k-1}g_{(\,\,)}f_{(i_2-1,\dots,i_k-1,n)}$
and
$g_{(0,i_2,\dots,i_k)}f_{(\,\,)}t_n=(-1)^{k-1}g_{(i_2-1,\dots,i_k-1,n)}f_{(\,\,)}$.
Given any fixed permutation $\sigma\in\Sigma_k$, consider the
summand $$(-1)^{{\rm sign}(\sigma)}
g_{(\widehat{\sigma(0)},\widehat{\sigma(i_2)},\dots,\widehat{\sigma(i_m)})}
f_{(\widehat{\sigma(i_{m+1})},\dots,\widehat{\sigma(i_k)})}t_n,\quad
m\geqslant 1,$$ on the right-hand side of $(1.6)$. Suppose that
$\sigma(0)>0$. In this case, we have
$\sigma(i_{m+1})=\widehat{\sigma(i_{m+1})}=0$. Therefore, taking
into account the relations $\widehat{\sigma(i_s)}=\sigma(i_s)$,
$m+2\leqslant s\leqslant k$, we obtain $$(-1)^{{\rm sign}(\sigma)}
g_{(\widehat{\sigma(0)},\widehat{\sigma(i_2)},\dots,\widehat{\sigma(i_m)})}
f_{(\widehat{\sigma(i_{m+1})},\dots,\widehat{\sigma(i_k)})}t_n=$$
$$=(-1)^{{\rm sign}(\sigma)+k-m-1}
g_{(\widehat{\sigma(0)},\widehat{\sigma(i_2)},\dots,\widehat{\sigma(i_m)})}
f_{(\sigma(i_{m+2})-1,\dots,\sigma(i_k)-1,n)}.$$ Let
$\varrho\in\Sigma_k$ be the permutation of the $k$-tuple
$(i_2-1,\dots,i_k-1,n)$ defined by
$$\varrho(i_2-1)=\sigma(0)-1,\quad
\varrho(i_3-1)=\sigma(i_2)-1,\dots,\varrho(i_{m+1}-1)=\sigma(i_m)-1,$$
$$\varrho(i_{m+2}-1)=\sigma(i_{m+2})-1,\dots,\varrho(i_k-1)=\sigma(i_k)-1,\quad
\varrho(n)=n.$$ Comparing the tuples
$(\sigma(0),\sigma(i_2),\dots,\sigma(i_k))$ and
$(\varrho(i_2-1),\dots,\varrho(i_k-1),\varrho(n))$, we see that
$$\widehat{\varrho(i_2-1)}=\widehat{\sigma(0)},\quad
\widehat{\varrho(i_3-1)}=\widehat{\sigma(i_2)},\dots,\widehat{\varrho(i_{m+1}-1)}=\widehat{\sigma(i_m)},$$
$$\widehat{\varrho(i_{m+2}-1)}=\sigma(i_{m+2})-1,\dots,\widehat{\varrho(i_k-1)}=\sigma(i_k)-1,\quad
\widehat{\varrho(n)}=n,$$ $${\rm sign}(\varrho)={\rm
sign}(\sigma)-m.$$ Since $\widehat{\varrho(i_2-1)}<
\dots<\widehat{\varrho(i_{m+1}-1)}$ and
$\widehat{\varrho(i_{m+2}-1)}<\dots<\widehat{\varrho(i_k-1)}<
\widehat{\varrho(n)}$, it follows that the right-hand side of
$(1.7)$ contains the summand $$(-1)^{{\rm sign}(\varrho)+k-1}
g_{(\widehat{\varrho(i_2-1)},\dots,\widehat{\varrho(i_{m+1}-1)})}
f_{(\widehat{\varrho(i_{m+2}-1)},\dots,\widehat{\varrho(i_k-1)},\widehat{\varrho(n)})}.$$
Clearly, this summands satisfies the relation $$(-1)^{{\rm
sign}(\varrho)+k-1}
g_{(\widehat{\varrho(i_2-1)},\dots,\widehat{\varrho(i_{m+1}-1)})}
f_{(\widehat{\varrho(i_{m+2}-1)},\dots,\widehat{\varrho(i_k-1)},\widehat{\varrho(n)})}=$$
$$=(-1)^{{\rm sign}(\sigma)}
g_{(\widehat{\sigma(0)},\widehat{\sigma(i_2)},\dots,\widehat{\sigma(i_m)})}
f_{(\widehat{\sigma(i_{m+1})},\dots,\widehat{\sigma(i_k)})}t_n.$$

Now, suppose that $\sigma(0)=0$. In this case, we have
$\widehat{\sigma(0)}=\sigma(0)=0$. Therefore, we obtain
$$(-1)^{{\rm sign}(\sigma)}
g_{(\widehat{\sigma(0)},\widehat{\sigma(i_2)},\dots,\widehat{\sigma(i_m)})}
f_{(\widehat{\sigma(i_{m+1})},\dots,\widehat{\sigma(i_k)})}t_n=$$
$$=(-1)^{{\rm
sign}(\sigma)+m-1}g_{(\widehat{\sigma(i_2)}-1,\dots,\widehat{\sigma(i_m)}-1,n-(k-m))}
f_{(\widehat{\sigma(i_{m+1})}-1,\dots,\widehat{\sigma(i_k)}-1)}.$$
Let $\varrho\in\Sigma_k$ be the permutation of the $k$-tuple
$(i_2-1,\dots,i_k-1,n)$ defined by
$$\varrho(i_2-1)=\sigma(i_2)-1,\dots,\varrho(i_m-1)=\sigma(i_m)-1,\quad
\varrho(i_{m+1}-1)=n,$$
$$\varrho(i_{m+2}-1)=\sigma(i_{m+1})-1,\dots,\varrho(i_k-1)=\sigma(i_{k-1})-1,\quad
\varrho(n)=\sigma(i_k)-1.$$ Comparing the tuples
$(\sigma(0),\sigma(i_2),\dots,\sigma(i_k))$ and
$(\varrho(i_2-1),\dots,\varrho(i_k-1),\varrho(n))$, we see that
$$\widehat{\varrho(i_2-1)}=\widehat{\sigma(i_2)}-1,\dots,
\widehat{\varrho(i_m-1)}=\widehat{\sigma(i_m)}-1,\quad\widehat{\varrho(i_{m+1}-1)}=n-(k-m),$$
$$\widehat{\varrho(i_{m+2}-1)}=\widehat{\sigma(i_{m+1})}-1,\dots,\widehat{\varrho(i_k-1)}=\sigma(i_{k-1})-1,\quad
\widehat{\varrho(n)}=\widehat{\sigma(i_k)}-1,$$ $${\rm
sign}(\varrho)={\rm sign}(\sigma)+(k-m).$$ Since
$\widehat{\varrho(i_2-1)}< \dots<\widehat{\varrho(i_{m+1}-1)}$ and
$\widehat{\varrho(i_{m+2}-1)}<\dots<\widehat{\varrho(i_k-1)}<
\widehat{\varrho(n)}$, it follows that the right-hand side of
$(1.7)$ contains the summand $$(-1)^{{\rm sign}(\varrho)+k-1}
g_{(\widehat{\varrho(i_2-1)},\dots,\widehat{\varrho(i_{m+1}-1)})}
f_{(\widehat{\varrho(i_{m+2}-1)},\dots,\widehat{\varrho(i_k-1)},\widehat{\varrho(n)})}.$$
Clearly, this summand satisfies the relation $$(-1)^{{\rm
sign}(\varrho)+k-1}
g_{(\widehat{\varrho(i_2-1)},\dots,\widehat{\varrho(i_{m+1}-1)})}
f_{(\widehat{\varrho(i_{m+2}-1)},\dots,\widehat{\varrho(i_k-1)},\widehat{\varrho(n)})}=$$
$$=(-1)^{{\rm sign}(\sigma)}
g_{(\widehat{\sigma(0)},\widehat{\sigma(i_2)},\dots,\widehat{\sigma(i_m)})}
f_{(\widehat{\sigma(i_{m+1})},\dots,\widehat{\sigma(i_k)})}t_n.$$
Thus, we have shown that each summand on the right-hand side of
$(1.6)$ is equal to a summand on the right-hand side of $(1.7)$.
It follows that the right-hand sides of $(1.6)$ and $(1.7)$ are
equal, because the number of summands on the right-hand side of
$(1.6)$ equals that on the right-hand side of $(1.7)$ and,
moreover, the permutations $\sigma$ and $\varrho$ uniquely
determine one another.~~~$\blacksquare$

It is clear that the associativity of the composition operation of
$\F$-modules implies the associativity of the composition
operation of $C\F$-modules. Moreover, for each $C\F$-module
$(X,d,\p,t)$, we have the identity morphism
$$1_X=\{(1_X)_{(i_1,\dots,i_k)}\}:(X,d,\p,t)\to (X,d,\p,t),$$
where $(1_X)_{(\,\,)}={\rm id_X}$ and $(1_X)_{(i_1,\dots,i_k)}=0$
for all $k\geqslant 1$. Thus, the class of all $C\F$-mo\-du\-les
over any commutative unital ring $K$ and their morphisms is a
category, which we denote $C\F(K)$.

Now, we recall that a differential homotopy or, more briefly, a
homotopy between morphisms $f,g:(X,d,\p_i)\to(Y,d,\p_i)$ of
differential modules with simplicial faces is defines as a
differential homotopy $h:X_{*,\bu}\to Y_{*,\bu+1}$ between
morphisms of differential modules $f,g:(X,d)\to (Y,d)$, which
satisfies the relations $\p_ih+h\p_i=0$, $i\geqslant 0$.

Let us consider the notion of a homotopy between morphisms of
differential modules with $\infty$-simplicial faces \cite{Lap1}
(see also \cite{Lap4}), which homotopically generalizes the notion
of a homotopy between morphisms of differential modules with
simplicial faces.

A homotopy between morphisms of $\F$-modules $f,g:(X,d,\p)\to
(Y,d,\p)$ is defined as a family of module maps
$h=\{h_{(i_1,\dots,i_k)}:X_{n,\bu}\to Y_{n-k,\bu+k+1}\}$,
$0\leqslant k\leqslant n$, $i_1,\dots,i_k\in\mathbb{Z}$,
$0\leqslant i_1<\dots<i_k\leqslant n$ (at $k=0$ we will use the
denotation $h_{(\,\,)}$), which satisfy the relations
$$d(h_{(i_1,\dots,i_k)})=f_{(i_1,\dots,i_k)}-g_{(i_1,\dots,i_k)}-
\partial_{(i_1,\dots,i_k)}h_{(\,\,)}-h_{(\,\,)}\partial_{(i_1,\dots,i_k)}\,+$$
$$+\sum_{\sigma\in\Sigma_k}\sum_{I_{\sigma}}(-1)^{{\rm
sign}(\sigma)+1}
\partial_{(\widehat{\sigma(i_1)},\dots,\widehat{\sigma(i_m)})}
h_{(\widehat{\sigma(i_{m+1})},\dots,\widehat{\sigma(i_k)})}\,+$$
$$+\,h_{(\widehat{\sigma(i_1)},\dots,\widehat{\sigma(i_m)})}
\partial_{(\widehat{\sigma(i_{m+1})},\dots,\widehat{\sigma(i_k)})},\eqno(1.8)$$
where $I_{\sigma}$ is the same as in $(1.1)$. The maps
$h_{(i_1,\dots ,i_k)}\in h$ are called the components of the
homotopy $h$.

For example, at $k=0,1,2,3$ the relations $(1.8)$ take,
respectively, the following view
$$d(h_{(\,\,)})=f_{(\,\,)}-g_{(\,\,)},\quad
d(h_{(i)})=f_{(i)}-g_{(i)}-\p_{(i)}h_{(\,\,)}-h_{(\,\,)}\p_{(i)},\quad
i\geqslant 0,$$ $$d(h_{(i,j)})=f_{(i,j)}-g_{(i,j)}-
\p_{(i,j)}h_{(\,\,)}-h_{(\,\,)}\p_{(i,j)}-\p_{(i)}h_{(j)}+\p_{(j-1)}h_{(i)}-
h_{(i)}\p_{(j)}+h_{(j-1)}\p_{(i)},~i<j,$$
$$d(h_{(i_1,i_2,i_3)})=f_{(i_1,i_2,i_3)}-g_{(i_1,i_2,i_3)}-\p_{(i_1,i_2,i_3)}h_{(\,\,)}-h_{(\,\,)}\p_{(i_1,i_2,i_3)}-
\p_{(i_1)}f_{(i_2,i_3)}-\p_{(i_1,i_2)}f_{(i_3)}\,-$$
$$-\,\p_{(i_3-2)}f_{(i_1,i_2)}-\p_{(i_2-1,i_3-1)}f_{(i_1)}+\p_{(i_2-1)}f_{(i_1,i_3)}+\p_{(i_1,i_3-1)}f_{(i_2)}-
h_{(i_1)}\p_{(i_2,i_3)}-h_{(i_1,i_2)}\p_{(i_3)}\,-$$
$$-\,h_{(i_3-2)}\p_{(i_1,i_2)}-
h_{(i_2-1,i_3-1)}\p_{(i_1)}+h_{(i_2-1)}\p_{(i_1,i_3)}+h_{(i_1,i_3-1)}\p_{(i_2)},\quad
i_1<i_2<i_3.$$

{\bf Definition 1.2}. A homotopy between an arbitrary morphisms of
$C\F$-modules $f,g:(X,d,\p,t)\to (Y,d,\p,t)$ is defined as any
homotopy $h=\{h_{(i_1,\dots,i_k)}\}$ between morphisms of
$\F$-modules $f,g:(X,d,\p)\to (Y,d,\p)$ whose components satisfy
the following conditions: $$h_{(\,\,)}t_n=t_nh_{(\,\,)},\quad
h_{(i_1,\dots,i_k)}t_n=\left\{\begin{array}{lll}
t_{n-k}h_{(i_1-1,\dots,i_k-1)},&k\geqslant 1,& i_1>0,\\
(-1)^{k-1}h_{(i_2-1,\dots,i_k-1,n)},&k\geqslant 1,& i_1=0.\\
\end{array}\right.\eqno(1.9)$$

{\bf Proposition 1.1}. For any $C\F$-modules $(X,d,\p,t)$ and
$(Y,d,\p,t)$, the relation between morphisms of $C\F$-modules of
the form $(X,d,\p,t)\to (Y,d,\p,t)$ defined by the presence of a
homotopy between them is an equivalence relation.

{\bf Proof}. Suppose given any morphism of $C\F$-modules
$f:(X,d,\p,t)\to (Y,d,\p,t)$. Then we have the homotopy
$0=\{0_{(i_1,\dots,i_k)}=0\}$ between morphisms $f$ and $f$.
Suppose given a homotopy $h=\{h_{(i_1,\dots,i_k)}\}$ between
morphisms $f:(X,d,\p,t)\to (Y,d,\p,t)$ and $g:(X,d,\p,t)\to
(Y,d,\p,t)$. Then the family of maps $-h=\{-h_{(i_1,\dots,i_k)}\}$
is a homotopy between morphisms $g$ and $f$. Suppose given a
homotopy $h=\{h_{(i_1,\dots,i_k)}\}$ between morphisms
$f:(X,d,\p,t)\to (Y,d,\p,t)$ and $g:(X,d,\p,t)\to (Y,d,\p,t)$ and,
moreover, given a homotopy $H=\{H_{(i_1,\dots,i_k)}\}$ between
morphisms $g:(X,d,\p,t)\to (Y,d,\p,t)$ and $p:(X,d,\p,t)\to
(Y,d,\p,t)$. Then the family of maps
$h+H=\{h_{(i_1,\dots,i_k)}+H_{(i_1,\dots,i_k)}\}$ is a homotopy
between morphisms $f$ and $p$.~~~$\blacksquare$

By using specified in Proposition 1.1 the equivalence relation
between morphisms of $C\F$-modules the notion of a homotopy
equivalence of $C\F$-modules is introduced in the usual way.
Namely, a morphism of $C\F$-modules is called a homotopy
equivalence of $C\F$-modules, when this morphism have a homotopy
inverse morphism of $C\F$-mo\-dules. \vspace{0.5cm}

\centerline{\bf \S\,2. The homotopy invariance of cyclic homology
of $C\F$-modules.}
\vspace{0.5cm}

First, recall that a $\D$-differential module \cite{Lap9} (sees
also \cite{Lap10}-\cite{Lap17}) or, more briefly, a $\D$-module
$(X,d^{\,i})$ is defined as a module $X$ together with a family of
module maps $\{d^{\,i}:X\to X~|~i\in\mathbb{Z},~i\geqslant 0\}$
satisfying the relations
$$\sum\limits_{i+j=k}d^{\,i}d^{\,j}=0,\quad k\geqslant
0.\eqno(2.1)$$

It is worth noting that a $\D$-module $(X,d^{\,i})$ can be
equipped with any $\mathbb{Z}^{\times n}$-gra\-ding, i.e.,
$X=\{X_{k_1,\dots,k_n}\}$, where
$(k_1,\dots,k_n)\in\mathbb{Z}^{\times n}$ and $n\geqslant 1$, and
the module maps $d^{\,i}:X\to X$ can have any $n$-degree
$(l_1(i),\dots,l_n(i))\in\mathbb{Z}^{\times n}$ for each
$i\geqslant 0$, i.e., $d^{\,i}:X_{k_1,\dots,k_n}\to
X_{k_1+l_1(i),\dots,k_n+l_n(i)}$.

For $k=0$, the relations $(2.1)$ have the form $d^{\,0}d^{\,0}=0$,
and hence $(X,d^{\,0})$ is a differential module. In \cite{Lap9}
the homotopy invariance of the $\D$-module structure over any
unital commutative ring under homotopy equivalences of
differential modules was established. Later, it was shown in
\cite{LV} that the homotopy invariance of the $\D$-mo\-du\-le
structure over fields of characteristic zero can be established by
using the Koszul duality theory.

It is also worth saying that in \cite{Lap9} by using specified
above homotopy invariance of the $\D$-differential module
structure the relationship between $\D$-differential modules and
spectral sequences was established. More precisely, in \cite{Lap9}
was shown that over an arbitrary field the category of
$\D$-differential modules is equivalent to the category of
spectral sequences.

Now, we recall \cite{Lap9} that a $\D$-module $(X,d^{\,i})$ is
said to be stable if, for any $x\in X$, there exists a number
$k=k(x)\geqslant 0$ such that $d^{\,i}(x)=0$ for each $i>k$. Any
stable $\D$-module $(X,d^{\,i})$ determines the differential
$\overline{d\,}:X\to X$ defined by
$\overline{d\,}=(d^{\,0}+d^{\,1}+\dots+d^{\,i}+\dots)$. The map
$\overline{d\,}:X\to X$ is indeed a differential because relations
$(2.1)$ imply the equality $\overline{d\,}\,\overline{d\,}=0$. It
is easy to see that if the stable $\D$-mo\-du\-le $(X,d^{\,i})$ is
equipped with a $\mathbb{Z}^{\times n}$-grading
$X=\{X_{k_1,\dots,k_n}\}$, where $k_1\geqslant
0,\dots,k_n\geqslant 0$, and maps $d^{\,i}:X\to X$, $i\geqslant
0$, have $n$-degree $(l_1(i),\dots,l_n(i))$ satisfying the
condition $l_1(i)+\dots+l_n(i)=-1$, then there is the chain
complex $(\overline{X},\overline{d\,})$ defined by the following
formulae:
$$\overline{X}_m=\bigoplus_{k_1+\dots+k_n=m}X_{k_1\dots,k_n},\quad
\overline{d\,}=\sum_{i=0}^\infty d^{\,i}:\overline{X}_m\to
\overline{X}_{m-1},\quad m\geqslant 0.$$

It was shown in \cite{Lap1} that any $\F$-module $(X,d,\p)$
determines the sequence of stable $\D$-modules $(X,d_q^{\,i})$,
$q\geqslant 0$, equipped with the bigrading $X=\{X_{n,m}\}$,
$n\geqslant 0$, $m\geqslant 0$, and defined by the following
formulae: $$d_q^{\,0}=d,~~ d_q^{\,k}=\sum\limits_{0\leqslant
i_1<\dots<i_k\leqslant n-q
}(-1)^{i_1+\dots+i_k}\p_{(i_1,\dots,i_k)}:X_{n,\bu}\to
X_{n-k,\bu+k-1},~~k\geqslant 1.\eqno(2.2)$$

Let us recall \cite{Lapin} the construction of the chain bicomplex
$(C(\overline{X}),\delta_1,\delta_2)$ that is defined by the
$C\F$-module $(X,d,\p,t)$. Given any $C\F$-mo\-dule $(X,d,\p,t)$,
consider the two $\D$-mo\-dules $(X,d_0^{\,i})$ and
$(X,d_1^{\,i})$ defined by $(2.2)$ for $q=0,1$, and the two
families of maps $$T_n=(-1)^nt_n:X_{n,\bu}\to X_{n,\bu},\quad
n\geqslant 0,$$ $$N_n=1+T_n+T_n^2+\dots+T_n^n:X_{n,\bu}\to
X_{n,\bu},\quad n\geqslant 0.$$ Obviously, the condition
$t_n^{\,n+1}=1$, $n\geqslant 0$, implies the relations
$$(1-T_n)N_n=0,\quad N_n(1-T_n)=0,\quad n\geqslant 0.\eqno(2.3)$$
Moreover, in \cite{Lapin} it was shown that the families of module
maps $\{T_n:X_{n,\bu}\to X_{n,\bu}\}$, $\{N_n:X_{n,\bu}\to
X_{n,\bu}\}$, $\{d_0^{\,i}:X_{*,\bu}\to X_{*-i,\bu+i-1}\}$ and
$\{d_1^{\,i}:X_{*,\bu}\to X_{*-i,\bu+i-1}\}$ are related by
$$d_0^{\,i}(1-T_n)=(1-T_{n-i})d_1^{\,i},\quad
d_1^{\,i}N_n=N_{n-i}d_0^{\,i},\quad i\geqslant 0,\quad n\geqslant
0.\eqno(2.4)$$

Now, we consider the chain complexes $(\overline{X},b)$ and
$(\overline{X},b^{'})$ corresponding to the $\D$-mo\-dules
$(X,d_0^{\,i})$ and $(X,d_1^{\,i})$ specified above. It is easy to
see that the chain complexes $(\overline{X},b)$ and
$(\overline{X},b^{'})$ are defined by
$$\overline{X}_n=\bigoplus_{k=0}^n X_{k,n-k},\quad
b=\overline{d\,}_0=\sum_{i=0}^n d_0^{\,i}:\overline{X}_n\to
\overline{X}_{n-1},$$ $$b^{'}=\overline{d\,}_1=\sum_{i=0}^n
d_1^{\,i}:\overline{X}_n\to \overline{X}_{n-1},\quad n\geqslant
0.$$ Consider also the two families of maps
$$\overline{T}_n=\sum_{k=0}^nT_{k}:\overline{X}_n\to\overline{X}_n,\quad
\overline{N}_n=\sum_{k=0}^nN_{k}:\overline{X}_n\to\overline{X}_n,\quad
n\geqslant 0.$$ The formulae $(2.3)$ and $(2.4)$ implies the
relations $$(1-\overline{T}_n)\overline{N}_n=0,\quad
\overline{N}_n(1-\overline{T}_n)=0,\quad n\geqslant 0,$$
$$b(1-\overline{T}_n)=(1-\overline{T}_{n-1})b^{'},\quad
b^{'}\overline{N}_n=\overline{N}_{n-1}b,\quad n\geqslant 0.$$ It
follows from these relations that any $C\F$-module $(X,d,\p,t)$
determines the chain bicomplex \vspace{1cm}

\begin{center}
\parbox {4cm}
{\setlength{\unitlength}{1cm}
\begin{picture}(2,3.5)

\put(-3.1,3.5){\makebox(0,0){\vspace{0.3cm}$\vdots$}}
\put(-0.7,3.5){\makebox(0,0){\vspace{0.3cm}$\vdots$}}
\put(1.7,3.5){\makebox(0,0){\vspace{0.3cm}$\vdots$}}
\put(4.1,3.5){\makebox(0,0){\vspace{0.3cm}$\vdots$}}

\put(-3.1,3.2){\vector(0,-1){0.7}}
\put(-0.7,3.2){\vector(0,-1){0.7}}
\put(1.7,3.2){\vector(0,-1){0.7}}
\put(4.1,3.2){\vector(0,-1){0.7}}

\put(-3,2){\makebox(0,0){$\overline{X}_{n+1}$}}
\put(-0.6,2){\makebox(0,0){$\overline{X}_{n+1}$}}
\put(1.8,2){\makebox(0,0){$\overline{X}_{n+1}$}}
\put(4.2,2){\makebox(0,0){$\overline{X}_{n+1}$}}
\put(6.3,2){\makebox(0,0){$\dots$}}

\put(-1.3,2){\vector(-1,0){1.1}} \put(1.1,2){\vector(-1,0){1.1}}
\put(3.5,2){\vector(-1,0){1.1}} \put(5.9,2){\vector(-1,0){1.1}}

\put(-3.3,2.8){\makebox(0,0){$^{b}$}}
\put(-3.3,1.1){\makebox(0,0){$^{b}$}}
\put(-3.3,-0.6){\makebox(0,0){$^{b}$}}
\put(-3.3,-2.3){\makebox(0,0){$^{b}$}}

\put(-1.0,2.8){\makebox(0,0){$^{-b^{'}}$}}
\put(-1.0,1.1){\makebox(0,0){$^{-b^{'}}$}}
\put(-1.0,-0.6){\makebox(0,0){$^{-b^{'}}$}}
\put(-1.0,-2.3){\makebox(0,0){$^{-b^{'}}$}}

\put(1.5,2.8){\makebox(0,0){$^{b}$}}
\put(1.5,1.1){\makebox(0,0){$^{b}$}}
\put(1.5,-0.6){\makebox(0,0){$^{b}$}}
\put(1.5,-2.3){\makebox(0,0){$^{b}$}}

\put(3.8,2.8){\makebox(0,0){$^{-b^{'}}$}}
\put(3.8,1.1){\makebox(0,0){$^{-b^{'}}$}}
\put(3.8,-0.6){\makebox(0,0){$^{-b^{'}}$}}
\put(3.8,-2.3){\makebox(0,0){$^{-b^{'}}$}}

\put(-3.1,1.7){\vector(0,-1){1.0}}
\put(-0.7,1.7){\vector(0,-1){1.0}}
\put(1.7,1.7){\vector(0,-1){1.0}}
\put(4.1,1.7){\vector(0,-1){1.0}}

\put(-3,0.2){\makebox(0,0){$\overline{X}_n$}}
\put(-0.6,0.2){\makebox(0,0){$\overline{X}_n$}}
\put(1.8,0.2){\makebox(0,0){$\overline{X}_n$}}
\put(4.2,0.2){\makebox(0,0){$\overline{X}_n$}}
\put(6.3,0.2){\makebox(0,0){$\dots$}}

\put(-1.3,0.2){\vector(-1,0){1.1}}
\put(1.1,0.2){\vector(-1,0){1.1}}
\put(3.5,0.2){\vector(-1,0){1.1}}
\put(5.9,0.2){\vector(-1,0){1.1}}

\put(-3.1,-0.1){\vector(0,-1){1.0}}
\put(-0.7,-0.1){\vector(0,-1){1.0}}
\put(1.7,-0.1){\vector(0,-1){1.0}}
\put(4.1,-0.1){\vector(0,-1){1.0}}

\put(-3,-1.55){\makebox(0,0){$\overline{X}_{n-1}$}}
\put(-0.6,-1.55){\makebox(0,0){$\overline{X}_{n-1}$}}
\put(1.8,-1.55){\makebox(0,0){$\overline{X}_{n-1}$}}
\put(4.2,-1.55){\makebox(0,0){$\overline{X}_{n-1}$}}
\put(6.3,-1.55){\makebox(0,0){$\dots$}}

\put(-1.7,-1.3){\makebox(0,0){$^{1-\overline{T}_{n-1}}$}}
\put(-1.7,0.4){\makebox(0,0){$^{1-\overline{T}_n}$}}
\put(-1.7,2.2){\makebox(0,0){$^{1-\overline{T}_{n+1}}$}}

\put(0.6,-1.3){\makebox(0,0){$^{\overline{N}_{n-1}}$}}
\put(0.6,0.4){\makebox(0,0){$^{\overline{N}_n}$}}
\put(0.6,2.2){\makebox(0,0){$^{\overline{N}_{n+1}}$}}

\put(3.1,-1.3){\makebox(0,0){$^{1-\overline{T}_{n-1}}$}}
\put(3.1,0.4){\makebox(0,0){$^{1-\overline{T}_n}$}}
\put(3.1,2.2){\makebox(0,0){$^{1-\overline{T}_{n+1}}$}}

\put(5.4,-1.3){\makebox(0,0){$^{\overline{N}_{n-1}}$}}
\put(5.4,0.4){\makebox(0,0){$^{\overline{N}_n}$}}
\put(5.4,2.2){\makebox(0,0){$^{\overline{N}_{n+1}}$}}

\put(-3.1,-2.8){\makebox(0,0){$\vdots$}}
\put(-0.7,-2.8){\makebox(0,0){$\vdots$}}
\put(1.7,-2.8){\makebox(0,0){$\vdots$}}
\put(4.1,-2.8){\makebox(0,0){$\vdots$}}

\put(-1.3,-1.55){\vector(-1,0){1.1}}
\put(1.1,-1.55){\vector(-1,0){1.1}}
\put(3.5,-1.55){\vector(-1,0){1.1}}
\put(5.9,-1.55){\vector(-1,0){1.1}}

\put(-3.1,-1.9){\vector(0,-1){0.7}}
\put(-0.7,-1.9){\vector(0,-1){0.7}}
\put(1.7,-1.9){\vector(0,-1){0.7}}
\put(4.1,-1.9){\vector(0,-1){0.7}}
\end{picture}}
\end{center}
\vspace{3.5cm}

\noindent We denote this chain bicomplex by
$(C(\overline{X}),D_1,D_2)$, where
$C(\overline{X})_{n,m}=\overline{X}_n$, $n\geqslant 0$,
$m\geqslant 0$, $D_1:C(\overline{X})_{n,m}\to
C(\overline{X})_{n,m-1}$, $D_2:C(\overline{X})_{n,m}\to
C(\overline{X})_{n-1,m}$, $$D_1=\left\{\begin{array}{ll}
1-\overline{T}_n,&m\equiv 1\,{\rm mod}(2),\\
\overline{N}_n,&m\equiv 0\,{\rm mod}(2).\\
\end{array}\right.\quad D_2=\left\{\begin{array}{ll}
b,&m\equiv 0\,{\rm mod}(2),\\ -b^{'},&m\equiv 1\,{\rm mod}(2),\\
\end{array}\right.$$
The chain complex associated with the chain bicomplex
$(C(\overline{X}),D_1,D_2)$ we denote by $({\rm
Tot}(C(\overline{X})),D)$, where $D=D_1+D_2$.

Recall \cite{Lapin} that the cyclic homology $HC(X)$ of a
$C\F$-module $(X,d,\p,t)$ is defined as the homology of the chain
complex $({\rm Tot}(C(\overline{X})),D)$ associated with the chain
bicomplex $(C(\overline{X}),D_1,D_2)$.

Now, we investigate functorial and homotopy properties of the
cyclic homology of $C\F$-mo\-du\-les.

{\bf Theorem 2.1}. The cyclic homology of $C\F$-modules over any
commutative unital ring $K$ determines the functor $HC:C\F(K)\to
GrM(K)$ from the category of $C\F$-mo\-dules $C\F(K)$ to the
category of graded $K$-modules $GrM(K)$. This functor sends
homotopy equivalences of $C\F$-modules into isomorphisms of graded
modules.

{\bf Proof}. First, show that every morphism of $C\F$-modules
induces a map of the graded cyclic homology modules. Suppose given
an arbitrary morphism of $C\F$-mo\-du\-les
$f=\{f_{(i_1,\dots,i_k)}\}:(X,d,\p,t)\to (X,d,\p,t)$. Consider the
family of maps $$f_q^k=\sum\limits_{0\leqslant
i_1<\dots<i_k\leqslant n-q
}(-1)^{i_1+\dots+i_k}f_{(i_1,\dots,i_k)}:X_{n,\bu}\to
Y_{n-k,\bu+k},\quad k\geqslant 0,\quad q\geqslant 0.\eqno(2.5)$$
For the family of maps $\{f_q^k\}$, by using $(1.3)$ we obtain the
relations
$$\sum_{i+j=k}d_q^{\,i}f_q^j=\sum_{i+j=k}f_q^id_q^{\,j},\quad
k\geqslant 0,\quad q\geqslant 0,\eqno(2.6)$$ where $(X,d_q^{\,i})$
and $(Y,d_q^{\,i})$ are sequences of $\D$-modules respectively
defined by $(2.2)$ for $\F$-modules $(X,d,\p)$ and $(Y,d,\p)$.
Similar to the way it was made in cite{Lapin}, direct calculations
with using $(1.5)$ show that the families maps $\{f_0^k\}$ and
$\{f_1^k\}$ satisfy the relations
$$f_0^k(1-T_n)=(1-T_{n-k})f_1^k,\quad f_1^kN_n=N_{n-k}f_0^k,\quad
k\geqslant 0,\quad n\geqslant 0.\eqno(2.7)$$ For $q=0,1$, the
equality $(2.6)$ imply that maps of graded modules
$$\overline{f}_0=\sum_{k=0}^nf_0^k:\overline{X}_n\to\overline{Y}_n,\quad
\overline{f}_1=\sum_{k=0}^nf_1^k:\overline{X}_n\to\overline{Y}_n,\quad
n\geqslant 0,$$ are chain maps $\overline{f}_0:(\overline{X},b)\to
(\overline{Y},b)$, $\overline{f}_1:(\overline{X},b^{'})\to
(\overline{Y},b^{'})$. From $(2.7)$ it follows that the chain maps
$\overline{f}_0$ and $\overline{f}_1$ satisfy the relations
$$\overline{f}_0(1-\overline{T}_n)=(1-\overline{T}_n)\overline{f}_1,\quad
\overline{f}_1\overline{N}_n=\overline{N}_n\overline{f}_0,\quad
n\geqslant 0.\eqno(2.8)$$ For chain bicomplexes
$(C(\overline{X}),D_1,D_2)$ and $(C(\overline{Y}),D_1,D_2)$,
consider the map of bigraded modules
$C(f):C(\overline{X})_{n,m}\to C(\overline{Y})_{n,m}$, $n\geqslant
0$, $m\geqslant 0$, defined by the following rule:
$$C(f)=\left\{\begin{array}{ll} \overline{f}_0,&m\equiv 0\,{\rm
mod}(2),\\ \overline{f}_1,&m\equiv 1\,{\rm mod}(2).\\
\end{array}\right.$$
From $(2.8)$ it follows that the map of bigraded modules
$C(f):C(\overline{X})\to C(\overline{Y})$ is the map of chain
bicomlexes $C(f):(C(\overline{X}),D_1,D_2)\to
(C(\overline{Y}),D_1,D_2)$. If we proceed to the homology of
associated chain complexes, then we obtain the map of graded
homology modules $$H({\rm Tot}(C(f))):H(({\rm
Tot}(C(\overline{X})),D))\to H(({\rm Tot}(C(\overline{Y})),D)).$$
Thus, every morphism of $C\F$-modules $f:(X,d,\p,t)\to (X,d,\p,t)$
induces the map of cyclic homology graded modules $HC(f)=H({\rm
Tot}(C(f))):HC(X)\to HC(Y)$.

Now, we consider the composition $gf:(X,d,\p,t)\to (Z,d,\p,t)$ of
morphisms of $C\F$-modules $f:(X,d,\p,t)\to (Y,d,\p,t)$ and
$g:(Y,d,\p,t)\to (Z,d,\p,t)$. Let us show that there is the
equality $HC(gf)=HC(g)HC(f)$ of maps of graded modules. Indeed,
since the composition $gf$ is a morphism of $C\F$-mo\-du\-les, we
have the family of maps $(gf)_q^k:X_{n,\bu}\to Z_{n,\bu}$,
$k\geqslant 0$, $q\geqslant 0$, defined by $(2.5)$. By using
$(1.4)$ it is easy to check that these maps satisfy the relations
$$(gf)_q^k=\sum\limits_{0\leqslant i_1<\dots<i_k\leqslant n-q
}\,\sum_{\sigma\in\Sigma_k}\sum_{I'_{\sigma}} (-1)^{{\rm
sign}(\sigma)+i_1+\dots+i_k}
g_{(\widehat{\sigma(i_1)},\dots,\widehat{\sigma(i_m)})}
f_{(\widehat{\sigma(i_{m+1})},\dots,\widehat{\sigma(i_k)})},$$
where $I'_\sigma$ is the same as in $(1.4)$. It is clear that
$$i_1+\dots+i_k+{\rm sign}(\sigma)\equiv
\widehat{\sigma(i_1)}+\dots+\widehat{\sigma(i_k)}\,{\rm mod}(2).$$
Moreover, for an arbitrary collections of integers $0\leqslant
j_1<\dots<j_m\leqslant n-q$ and $0\leqslant
j_{m+1}<\dots<j_k\leqslant n-q$, there is the collection
$0\leqslant i_1<\dots<i_k\leqslant n-q$ and the permutation
$\sigma\in\Sigma_k$ such that $\widehat{\sigma(i_s)}=j_s$ for each
$1\leqslant s\leqslant k$. Therefore specified above relations
imply the relations $$(gf)_q^k=\sum_{i+j=k}g_q^if_q^j,\quad
k\geqslant 0,\quad q\geqslant 0.$$ For $q=0,1$, by using these
relations we obtain the equalities
$\overline{(gf)}_0=\overline{g}_0\overline{f}_0$ and
$\overline{(gf)}_1=\overline{g}_1\overline{f}_1$. These equalities
imply the equality $C(gf)=C(g)C(f)$ of maps of chain bicomplexes.
If we proceed to the homology of associated chain complexes, then
we obtain the required equality $HC(gf)=HC(g)HC(f)$ of maps of
graded cyclic homology modules.

Suppose given an arbitrary homotopy $h$ between given morphisms of
$C\F$-modules $f:(X,d,t,\p)\to (Y,d,t,\p)$ and $g:(X,d,t,\p)\to
(Y,d,t,\p)$. In the same manner as above we obtain the homotopy
$C(h):C(\overline{X})_{*,\bu}\to C(\overline{Y})_{*+1,\bu}$
between maps of chain bicomplexes $C(f)$ and $C(g)$. If we proceed
to the homology of associated chain complexes, then we obtain the
equality $HC(f)=HC(g):HC(X)\to HC(Y)$ of maps of graded cyclic
homology modules. This equality implies that if the morphism of
$C\F$-modules $f:(X,d,\p,t)\to (Y,d,\p,t)$ is a homotopy
equivalence of $C\F$-mo\-dules, then the induces map
$HC(f):HC(X)\to HC(Y)$ of graded cyclic homology modules is an
isomorphism of graded modules.~~~$\blacksquare$
\vspace{0.5cm}

\centerline{\bf \S\,3. The homotopy invariance of cyclic homology
of $\A$-algebras.}
\vspace{0.5cm}

First, following \cite{Kad} and \cite{Smir} (see also \cite{S}),
we recall necessary definitions related to the notion of an
$\A$-algebra.

An $\A$-algebra $(A,d,\pi_n)$ is any differential module $(A,d)$
with $A=\{A_n\}$, $n\in\mathbb{Z}$, $n\geqslant 0$, $d:A_\bu\to
A_{\bu-1}$, equipped with a family of maps $\{\pi_n:(A^{\otimes
(n+2)})_\bu\to A_{\bu+n}\}$, $n\geqslant 0$, satisfying the
following relations for any integer $n\geqslant -1$:
$$d(\pi_{n+1})=\sum\limits_{m=0}^n\sum_{t=1}^{m+2}(-1)^{t(n-m+1)+n+1}
\pi_m(\underbrace{1\otimes\dots\otimes
1}_{t-1}\otimes\,\pi_{n-m}\otimes\underbrace{1\otimes\dots \otimes
1}_{m-t+2}),\eqno(3.1)$$ where
$d(\pi_{n+1})=d\pi_{n+1}+(-1)^n\pi_{n+1}d$. For example, at
$n=-1,0,1$ the relations $(3.1)$ take the forms $$d(\pi_0)=0,\quad
d(\pi_1)=\pi_0(\pi_0\otimes 1)-\pi_0(1\otimes\pi_0),$$
$$d(\pi_2)=\pi_0(\pi_1\otimes
1)+\pi_0(1\otimes\pi_1)-\pi_1(\pi_0\otimes 1\otimes
1)+\pi_1(1\otimes\pi_0\otimes 1)-\pi_1(1\otimes 1\otimes\pi_0).$$

A morphism of $\A$-algebras $f:(A,d,\pi_n)\to (A',d,\pi_n)$ is
defined as a family of module maps
$f=\{f_n:(A^{\otimes(n+1)})_\bu\to
A'_{\bu+n}~|~n\in\mathbb{Z},~n\geqslant 0\}$, which, for all
integers $n\geqslant -1$, satisfy the relations
$$d(f_{n+1})=\sum_{m=0}^n\sum_{t=1}^{m+1}(-1)^{t(n-m+1)+n+1}f_m(\underbrace{1\otimes\dots\otimes
1}_{t-1}\otimes\,\pi_{n-m}\otimes\underbrace{1\otimes\dots\otimes
1}_{m-t+1})\,-$$
$$-\sum_{m=0}^n\sum_{\,\,\,J_{n-m}}(-1)^{\varepsilon}\pi_m(f_{n_1}\otimes
f_{n_2}\otimes\dots\otimes f_{n_{m+2}}),\eqno(3.2)$$ where
$d(f_{n+1})=df_{n+1}+(-1)^nf_{n+1}d$ and $$J_{n-m}=\{n_1\geqslant
0,n_2\geqslant 0,\dots,n_{m+2}\geqslant
0~|~n_1+n_2+\dots+n_{m+2}=n-m\};$$
$$\varepsilon=\sum_{i=1}^{m+1}(n_i+1)(n_{i+1}+\dots+n_{m+2}).$$
For example, at $n=-1,0,1$ the relations $(3.2)$ take,
respectively, the following view $$d(f_0)=0,\quad
d(f_1)=f_0\pi_0-\pi_0(f_0\otimes f_0),$$
$$d(f_2)=f_0\pi_1-f_1(\pi_0\otimes
1)+f_1(1\otimes\pi_0)-\pi_0(f_1\otimes f_0)+\pi_0(f_0\otimes
f_1)-\pi_1(f_0\otimes f_0\otimes f_0).$$

Under a composition of morphisms of $\A$-algebras
$f:(A,d,\pi_n)\to(A',d,\pi_n)$ and
$g:(A',d,\pi_n)\to(A'',d,\pi_n)$ we mean the morphism of
$\A$-algebras $$gf=\{(gf)_n\}:(A,d,\pi_n)\to(A'',d,\pi_n)$$
defined by
$$(gf)_{n+1}=\sum_{m=-1}^n\sum_{J_{n-m}}(-1)^\varepsilon
g_{m+1}(f_{n_1}\otimes f_{n_2}\otimes\dots\otimes
f_{n_{m+2}}),\quad n\geqslant -1,\eqno(3.3)$$ where $J_{n-m}$ and
${\varepsilon}$ are the same as in $(3.2)$. For example, at
$n=0,1,2$ the formulae $(3.3)$ take, respectively, the following
view $$(gf)_0=g_0f_0,\quad (gf)_1=g_0f_1+g_1(f_0\otimes f_0),$$
$$(gf)_2=g_0f_2-g_1(f_0\otimes f_1)+g_1(f_1\otimes
f_0)+g_2(f_0\otimes f_0\otimes f_0).$$

It is easy to see that a composition of morphisms of $\A$-algebras
is associative. Moreover, for any $\A$-algebra $(A,d,\pi_n)$,
there is the identity morphism $$1_A=\{(1_A)_n\}:(A,d,\pi_n)\to
(A,d,\pi_n),$$ where $(1_A)_0={\rm id_A}$ and $(1_A)_n=0$ for each
$n\geqslant 1$. Thus, the class of all $\A$-algebras over any
commutative unital ring $K$ and their morphisms is a category,
which we denote $\A(K)$.

A homotopy between morphisms of $\A$-algebras $f,g:(A,d,\pi_n)\to
(A',d,\pi_n)$ is defined as a family of module maps
$h=\{h_n:(A^{\otimes(n+1)})_\bu\to
A'_{\bu+n+1}~|~n\in\mathbb{Z},~n\geqslant 0\}$, which, for all
integers $n\geqslant -1$, satisfy the relations
$$d(h_{n+1})=f_{n+1}-g_{n+1}+\sum_{m=0}^n\sum_{t=1}^{m+1}(-1)^{t(n-m+1)+n}h_m(\underbrace{1\otimes\dots\otimes
1}_{t-1}\otimes\,\pi_{n-m}\otimes\underbrace{1\otimes\dots\otimes
1}_{m-t+1})\,+$$
$$+\sum_{m=0}^n\sum_{\,\,\,J_{n-m}}\sum_{i=1}^{m+2}(-1)^\varrho\pi_m(g_{n_1}\otimes\dots\otimes
g_{n_{i-1}}\otimes h_{n_i}\otimes f_{n_{i+1}}\otimes\dots\otimes
f_{n_{m+2}}),\eqno(3.4)$$ where
$d(h_{n+1})=dh_{n+1}+(-1)^{n+1}h_{n+1}d$ and $J_{n-m}$ is the same
as in $(3.2)$;
$$\varrho=m+\sum_{k=1}^{m+1}(n_k+1)(n_{k+1}+\dots+n_{m+2})+\sum_{k=1}^{i-1}n_k.$$

For example, at $n=-1,0,1$ the relations $(3.4)$ take,
respectively, the following view $$d(h_0)=f_0-g_0,\quad
d(h_1)=f_1-g_1-h_0\pi_0+\pi_0(h_0\otimes f_0)+\pi_0(g_0\otimes
f_0),$$ $$d(h_2)=f_2-g_2-h_0\pi_1+h_1(\pi_0\otimes
1)-h_1(1\otimes\pi_0)-\pi_0(h_0\otimes f_1)-\pi_0(g_0\otimes
h_1)\,+$$ $$+\,\pi_0(h_1\otimes f_0)-\pi_0(g_1\otimes
h_0)-\pi_1(h_0\otimes f_0\otimes f_0)-\pi_1(g_0\otimes h_0\otimes
f_0)-\pi_1(g_0\otimes g_0\otimes h_0).$$

The origin of the signs in formulae $(3.1)$\,--\,$(3.4)$ is
described in detail in \cite{Lapin2}.

For any $\A$-algebras $(A,d,\pi_n)$ and $(A',d,\pi_n)$, the
relation between morphisms of $\A$-algebras of the form
$(A,d,\pi_n)\to (A',d,\pi_n)$ defined by the presence of a
homotopy between them is an equivalence relation. By using this
equivalence relation between morphisms of $\A$-algebras the notion
of a homotopy equivalence of $\A$-algebras is introduced in the
usual way. Namely, a morphism of $\A$-algebras is called a
homotopy equivalence of $\A$-algebras, when this morphism have a
homotopy inverse morphism of $\A$-algebras.

Now, let us proceed to cyclic homology of $\A$-algebras. In
\cite{Lapin} it was shown that any $\A$-algebra defines the tensor
$C\F$-module $({\cal L}(A),d,\p,t)$, which given by the following
equalities: $${\cal L}(A)=\{{\cal L}(A)_{n,m}\},\quad {\cal
L}(A)_{n,m}=(A^{\otimes (n+1)})_m,\quad n\geqslant 0,\quad
m\geqslant 0,$$ $$d(a_0\otimes\dots\otimes
a_n)=\sum_{i=0}^n(-1)^{|a_0|+\dots+|a_{i-1}|}a_0\otimes\dots\otimes
a_{i-1}\otimes d(a_i)\otimes a_{i+1}\otimes\dots\otimes a_n,$$
$$t_n(a_0\otimes\dots\otimes
a_n)=(-1)^{|a_n|(|a_0|+\dots+|a_{n-1}|)}a_n\otimes
a_0\otimes\dots\otimes a_{n-1},$$ where $|a|=q$ means that $a\in
A_q$. The family of module maps $$\p=\{\p_{(i_1,\dots,i_k)}:{\cal
L}(A)_{n,p}\to {\cal L}(A)_{n-k,p+k-1}\},$$ $$n\geqslant 0,\quad
p\geqslant 0,\quad 1\leqslant k\leqslant n,\quad 0\leqslant
i_1<\dots<i_k\leqslant n,$$ is defined by
$$\p_{(i_1,\dots,i_k)}\!=$$ $$=\left\{\begin{array}{ll}
(-1)^{k(p-1)}1^{\otimes j}\otimes\,\pi_{k-1}\otimes 1^{\otimes
(n-k-j)}\,,\,\,\mbox{if}\,\,0\leqslant j\leqslant
n-k&\\~~~\mbox{and}\,\,\, (i_1,\dots,i_k)=(j,j+1,\dots,j+k-1);\\
(-1)^{q(k-1)}\p_{(0,1,\dots,k-1)}t_n^q\,,\,\,\mbox{if}\,\,1\leqslant
q\leqslant k&\\~~~\mbox{and}\,\,\,
(i_1,\dots,i_k)=(0,1,\dots,k-q-1,n-q+1,n-q+2,\dots,n);\\
0,\quad\mbox{otherwise}.&\\
\end{array}\right.\eqno(3.5)$$

Recall \cite{Lapin} that the cyclic homology $HC(A)$ of an
$\A$-algebra $(A,d,\pi_n)$ is defined as the cyclic homology
$HC({\cal L}(A))$ of the $C\F$-module $({\cal L}(A),d,\p,t)$.

Now, we investigate functorial and homotopy properties of the
cyclic homology of $\A$-al\-geb\-ras.

{\bf Theorem 3.1}. The cyclic homology of $\A$-algebras over an
arbitrary commutative unital ring $K$ determines the functor
$HC:\A(K)\to GrM(K)$ from the category of $\A$-al\-gebras $\A(K)$
to the category of graded $K$-modules $GrM(K)$. This functor sends
homotopy equivalences of $\A$-algebras into isomorphisms of graded
modules.

{\bf Proof}. First, show that every morphism of $\A$-algebras
induces a morphism of $C\F$-modules. Given any morphism of
$\A$-algebras $f:(A,d,\pi_n)\to (A',d,\pi_n)$, we define the
family of module maps $${\cal L}(f)=\{{\cal
L}(f)^n_{(i_1,\dots,i_k)}:{\cal L}(A)_{n,p}\to {\cal
L}(A')_{n-k,p+k}\},$$ $$n\geqslant 0,\quad p\geqslant 0,\quad
0\leqslant k\leqslant n,\quad 0\leqslant i_1<\dots<i_k\leqslant
n,$$ by the following rules:

1). If $k=0$, then $${\cal L}(f)^n_{(\,\,)}=f_0^{\otimes
(n+1)}.\eqno(3.6)$$

2). If $i_k<n$ and
$(i_1,\dots,i_k)=((j^1_1,\dots,j^1_{n_1}),(j^2_1,\dots,j^2_{n_2}),\dots,(j^s_1,\dots,j^s_{n_s}))$,
$$1\leqslant s\leqslant k,\quad n_1\geqslant 1,\dots,n_s\geqslant
1,\quad n_1+\dots+n_s=k,$$ $$j_{b+1}^a=j_b^a+1,\quad 1\leqslant
a\leqslant s,\quad 1\leqslant b\leqslant n_a-1,\quad
j_1^{c+1}\geqslant j^c_{n_c}+2,\quad 1\leqslant c\leqslant s-1,$$
then $${\cal
L}(f)^n_{(i_1,\dots,i_k)}=(-1)^{k(p-1)+\gamma}\underbrace{f_0\otimes\dots\otimes
f_0}_{k_1}\otimes\,
f_{n_1}\!\otimes\underbrace{f_0\otimes\dots\otimes
f_0}_{k_2}\otimes \,f_{n_2}\,\otimes$$
$$\otimes\underbrace{f_0\otimes\dots\otimes
f_0}_{k_3}\otimes\,f_{n_3}\!\otimes\dots\otimes\underbrace{f_0\otimes\dots\otimes
f_0}_{k_s}\otimes\,f_{n_s}\!\otimes\underbrace{f_0\otimes\dots\otimes
f_0}_{k_{s+1}},\eqno(3.7)$$ where $k_1=j_1^1$,
$k_i=j^{\,i}_1-j^{\,i-1}_{n_{i-1}}-2$ at $2\leqslant i\leqslant
s$, $k_{s+1}=n+1-(k_1+\dots+k_s)-k-s$ and
$$\gamma=\sum\limits_{i=1}^{s-1}n_i(n_{i+1}+\dots+n_s).$$

3). If $i_k=n$ and
$$(i_1,\dots,i_k)=((0,1,\dots,z-1-q),(j^1_1-q,\dots,j^1_{n_1}-q),(j^2_1-q,\dots,j^2_{n_2}-q),\dots$$
$$\dots,(j^s_1-q,\dots,j^s_{n_s}-q),(n-q+1,n-q+2,\dots,n)),$$
$$z\geqslant 1,\quad  1\leqslant q\leqslant z,\quad 0\leqslant
s\leqslant k-1,\quad n_1\geqslant 1,\dots,n_s\geqslant 1,$$
$$z+n_1+\dots+n_s=k,\quad j_{b+1}^a=j_b^a+1,\quad 1\leqslant
a\leqslant s,\quad 1\leqslant b\leqslant n_a-1,$$ $$j_1^1\geqslant
z+1,\quad j_1^{c+1}\geqslant j^c_{n_c}+2,\quad 1\leqslant
c\leqslant s-1,\quad j^s_{n_s}\leqslant n-1,$$ then $${\cal
L}(f)^n_{(i_1,\dots,i_k)}=(-1)^{q(z-1)}{\cal
L}(f)^n_{((0,1,\dots,z-1),(j^1_1,\dots,j^1_{n_1}),\dots,(j^s_1,
\dots,j^s_{n_s}))}t_n^q.\eqno(3.8)$$

For example, if we consider the map ${\cal
L}(f)^{15}_{((2,3),(6,7,8))}:{\cal L}(A)_{15,p}\to {\cal
L}(A')_{10,p+5}$, then by $(3.7)$ we obtain $${\cal
L}(f)^{15}_{((2,3),(6,7,8))}=(-1)^{5(p-1)+2\cdot 3}f_0\otimes
f_0\otimes f_2\otimes f_0\otimes
f_3\otimes\underbrace{f_0\otimes\dots\otimes f_0}_6.$$ If we
consider consider the map ${\cal
L}(f)^n_{((0,1),(3,4),(n-2,n-1,n)}:{\cal L}(A)_{n,p}\to {\cal
L}(A')_{n-7,p+7}$, where $n\geqslant 8$, then by $(3.8)$ and
$(3.7)$ we obtain $${\cal L}(f)^n_{((0,1),(3,4),(n-2,n-1,n)}=$$
$$=(-1)^{3(5-1)}{\cal
L}(f)^n_{((0,1,2,3,4),(6,7))}t^3_n=(-1)^{7(p-1)+5\cdot
2}(f_5\otimes f_2\otimes\underbrace{f_0\otimes\dots\otimes
f_0}_{n-8})t^3_n.$$

It is worth noting that any collection of integers
$(i_1,\dots,i_k)$, $0\leqslant i_1<\dots<i_k\leqslant n$, always
can be written in the form specified in the rule 2) or in the rule
3).

Now, we show that the maps ${\cal L}(f)^n_{(i_1,\dots,i_k)}$
satisfy the conditions $(1.5)$. It is clear that at $k=0$ the the
equality ${\cal L}(f)^n_{(\,\,)}t_n=t_n{\cal L}(f)^n_{(\,\,)}$ is
true. Consider the maps ${\cal L}(f)^n_{(i_1,\dots,i_k)}$, where
$i_1>0$ and $i_k<n$, defined by formulae $(3.7)$ at $s=1$. Suppose
that
$(i_1,\dots,i_k)=(j^1_1,\dots,j^1_{n_1})=(j,j+1,\dots,j+k-1)$,
where $1\leqslant j\leqslant n-k$. In this case, on the one hand,
we have at any element $a_0\otimes\dots\otimes a_n\in {\cal
L}(A)_{n,p}=(A^{n+1})_p$ the equalities $${\cal
L}(f)^n_{(j,j+1,\dots,j+k-1)}t_n(a_0\otimes\dots\otimes
a_n)=(-1)^\alpha {\cal L}(f)^n_{(j,j+1,\dots,j+k-1)}(a_n\otimes
a_0\otimes\dots\otimes a_{n-1})=$$
$$=(-1)^{\alpha+k(p-1)}(f_{0}^{\otimes j}\otimes f_k\otimes
f_0^{\otimes (n-k-j)})(a_n\otimes a_0\otimes\dots\otimes
a_{n-1})=(-1)^{\alpha+k(p-1)+\beta}f_0(a_n)\,\otimes$$
$$\otimes\,f_0(a_0)\otimes\dots\otimes f_0(a_{j-2})\otimes
f_k(a_{j-1}\otimes\dots\otimes a_{j+k-1})\otimes
f_0(a_{j+k})\otimes\dots\otimes f_0(a_{n-1}),$$ where
$\alpha=|a_n|(|a_0|+\dots+|a_{n-1}|)$ and
$\beta=k(|a_n|+|a_0|+\dots+|a_{j-2}|)$. On the other hand, we have
the equalities $$t_{n-k}{\cal
L}(f)^n_{(j-1,j,\dots,j+k-2)}(a_0\otimes\dots\otimes
a_n)=(-1)^{k(p-1)}t_{n-k}(f_0^{\otimes (j-1)}\otimes f_k\otimes
f_0^{\otimes (n-k-j+1)})(a_0\,\otimes$$ $$\otimes\dots\otimes
a_n)=(-1)^{k(p-1)+\varphi}t_{n-k}(f_0(a_0)\otimes\dots\otimes
f_0(a_{j-2})\otimes f_k(a_{j-1}\otimes\dots\otimes
a_{j+k-1})\,\otimes$$ $$\otimes\,f_0(a_{j+k})\otimes\dots\otimes
f_0(a_n))=(-1)^{k(p-1)+\varphi+\delta} f_0(a_n)\otimes
f_0(a_0)\otimes\dots\otimes f_0(a_{j-2})\,\otimes$$
$$\otimes\,f_k(a_{j-1}\otimes\dots\otimes a_{j+k-1})\otimes
f_0(a_{j+k})\otimes\dots\otimes f_0(a_{n-1}),$$ where
$\varphi=k(|a_0|+\dots+|a_{j-2}|)$ and
$\delta=|a_n|(|a_0|+\dots+|a_{n-1}|+k)$. Since
$\alpha+\beta=\varphi+\delta$, we obtain the required relation
$${\cal L}(f)^n_{(j,j+1,\dots,j+k-1)}t_n=t_{n-k}{\cal
L}(f)^n_{(j-1,j,\dots,j+k-2)}.$$ In the similar way it is checked
that relations $(1.5)$ holds for all maps ${\cal
L}(f)^n_{(i_1,\dots,i_k)}$, where $i_1>0$ and $i_k<n$. Now,
consider the maps ${\cal L}(f)^n_{(i_1,\dots,i_k)}$, where $i_1=0$
and $i_k<n$, defined by $(3.7)$ at $s=1$. Suppose that
$(i_1,\dots,i_k)=(0,1,\dots,k-1)$, where $1\leqslant k\leqslant
n$. In this case, by using $(3.8)$ at $s=0$ and $q=1$ we obtain
the the required relation $${\cal
L}(f)^n_{(0,1,\dots,k-1)}t_n=(-1)^{k-1}{\cal
L}(f)^n_{(0,1,\dots,k-2,n)}.$$ In the similar way it is checked
that relations $(1.5)$ holds for all defined at $s\geqslant 2$ by
$(3.7)$ maps ${\cal L}(f)^n_{(i_1,\dots,i_k)}$, where $i_1=0$ and
$i_k<n$. Now, consider the maps ${\cal L}(f)^n_{(i_1,\dots,i_k)}$,
where $i_1=0$ and $i_k=n$, defined by $(3.8)$ at $s=0$. Suppose
that
$$(i_1,\dots,i_k)=((0,1,\dots,k-1-q),(n-q+1,n-q+2,\dots,n)),$$
where $1\leqslant q\leqslant k-1$. In this case, by using $(3.8)$
at $s=0$ we obtain $${\cal
L}(f)^n_{(0,1,\dots,k-q-1,n-q+1,n-q+2,\dots,n)}t_n=(-1)^{q(k-1)}{\cal
L}(f)^n_{(0,1,\dots,k-1)}t_n^{q+1}=$$ $$=(-1)^{k-1}{\cal
L}(f)^n_{(0,1,\dots,k-q-2,n-q,n-q+1,\dots,n-1,n)}.$$ In the
similar way it is proved that relations $(1.5)$ holds for all
defined at $s\geqslant 1$ by $(3.8)$ maps ${\cal
L}(f)^n_{(i_1,\dots,i_k)}$, where $i_1=0$ and $i_k=n$. Now,
consider the maps ${\cal L}(f)^n_{(i_1,\dots,i_k)}$, where $i_1>0$
and $i_k=n$, defined by $(3.8)$ at $s=0$. Suppose that
$(i_1,\dots,i_k)=(n-k+1,n-k+2,\dots,n)$. Then, by using $(3.8)$ at
$s=0$ and also applying the relations $t_{n-k}^{n-k+1}=1$ and
$t_n^{n+1}=1$, we obtain $${\cal
L}(f)^n_{(n-k+1,n-k+2,\dots,n)}t_n=(-1)^{k(k-1)}{\cal
L}(f)^n_{(0,1,\dots,k-1)}t_n^{k+1}= t_{n-k}^{n-k+1}{\cal
L}(f)^n_{(0,1,\dots,k-1)}t_n^{k+1}=$$ $$=t_{n-k}{\cal
L}(f)^n_{(n-k,n-k+1,\dots,n-1)}t_n^{n+1}=t_{n-k}{\cal
L}(f)^n_{(n-k,n-k+1,\dots,n-1)}.$$ In a similar manner is checked
that relations $(1.5)$ holds for all defined at $s\geqslant 1$ by
$(3.8)$ maps ${\cal L}(f)^n_{(i_1,\dots,i_k)}$, where $i_1>1$ and
$i_k=n$. Thus, all maps ${\cal L}(f)^n_{(i_1,\dots,i_k)}\in{\cal
L}(f)$ satisfy the relations $(1.5)$.

Now, let us show that the family of maps ${\cal L}(f)=\{{\cal
L}(f)^n_{(i_1,\dots,i_k)}\}$ is a morphism of $\F$-mo\-du\-les
${\cal L}(f):({\cal L}(A),d,\p)\to ({\cal L}(A'),d,\p)$. We must
check relations $(1.3)$ for the maps ${\cal
L}(f)^n_{(i_1,\dots,i_k)}\in{\cal L}(f)$. It is clear that at
$k=0$ we have $d({\cal L}(f)^n_{(\,\,)})=0$ because $d(f_0)=0$.
Now, we check that the maps $${\cal
L}(f)^{n+1}_{(0,1,\dots,n)}=(-1)^{(n+1)(p-1)}f_{n+1}:(A^{\otimes
(n+2)})_p\to A'_{p+n+1},\quad n\geqslant 0,$$ satisfy the
relations $(1.3)$. With this purpose we write the relations
$(3.2)$ in the form $$d(f_{n+1})=f_0\pi_n
+\sum_{m=0}^{n-1}\sum_{t=1}^{m+2}(-1)^{t(n-m)+n+1}f_{m+1}(\underbrace{1\otimes\dots\otimes
1}_{t-1}\otimes\,\pi_{n-m-1}\otimes\underbrace{1\otimes\dots\otimes
1}_{m+2-t})\,-$$ $$-\,\pi_n(f_0\otimes\dots\otimes f_0)-
\sum_{m=0}^{n-1}\sum_{s=1}^{m+2}\sum_{N_{n-m}}\sum_{T_{m+2}}(-1)^\mu\pi_m(\underbrace{f_0\otimes\dots\otimes
f_0}_{t_1-1}\otimes \,f_{n_1}\otimes$$
$$\otimes\underbrace{f_0\otimes\dots\otimes
f_0}_{t_2-1}\otimes\,f_{n_2}\otimes\dots\otimes\underbrace{f_0\otimes\dots\otimes
f_0}_{t_s-1}\otimes
\,f_{n_s}\!\otimes\underbrace{f_0\otimes\dots\otimes
f_0}_{m+2-(t_1+\dots+t_s)}),\quad n\geqslant -1,\eqno(3.9)$$ where
$$N_{n-m}=\{n_1\geqslant 1,n_2\geqslant 1,\dots,n_s\geqslant
1~|~n_1+n_2+\dots+n_s=n-m\},$$ $$T_{m+2}=\{t_1\geqslant
1,\dots,t_s\geqslant 1~|~t_1+\dots+t_s\leqslant m+2\},$$
$$\mu=\sum_{i=1}^s(t_i-1)(n_i+\dots+n_s)+\sum_{i=1}^{s-1}(n_i+1)(n_{i+1}+\dots+n_s).$$
Given any fixed collections $(n_1,\dots,n_s)\in I_{n-m}$ and
$(t_1,\dots,t_s)\in T_{m+2}$, consider a partition of the
collection $(0,1,\dots,n)$ into $2s+1$ blocks as
$$(0,1,\dots,n)=(a_1,b_1,a_2,b_2,\dots,a_s,b_s,a_{s+1}),$$
$$a_1=(0,1,\dots,t_1-2),\quad b_1=(t_1-1,t_1,\dots,t_1+n_1-2),$$
$$a_i=(\sum_{k=1}^{i-1}t_k+n_k-1,\sum_{k=1}^{i-1}t_k+n_k,\dots,\sum_{k=1}^{i-1}t_k+n_k+t_i-2),\quad
2\leqslant i\leqslant s,$$
$$b_i=(\sum_{k=1}^{i-1}t_k+n_k+t_i-1,\sum_{k=1}^{i-1}t_k+n_k+t_i,\dots,\sum_{k=1}^{i-1}t_k+n_k+t_i+n_i-2),\quad
2\leqslant i\leqslant s,$$
$$a_{s+1}=(\sum_{k=1}^st_k+n_k-1,\sum_{k=1}^st_k+n_k,\dots,n).$$
Given any specified above partition
$(0,1,\dots,n)=(a_1,b_1,a_2,b_2,\dots,a_s,b_s,a_{s+1})$, we define
the permutation $\sigma_{n_1,\dots
n_s,t_1,\dots,t_s}\in\Sigma_{n+1}$, which acting on the collection
of numbers $(0,1,\dots,n)$ by the following rule:
$$\sigma_{n_1,\dots
n_s,t_1,\dots,t_s}(0,1,\dots,n)=(a_1,a_2,\dots,a_s,a_{s+1},b_1,b_2,\dots,b_s).\eqno(3.10)$$
By using the relation $n_1+\dots+n_s=n-m$ it is easy verify that
the equality of collections $$(\widehat{\sigma_{n_1,\dots
n_s,t_1,\dots,t_s}(0)},\dots,\widehat{\sigma_{n_1,\dots
n_s,t_1,\dots,t_s}(n))}=(0,1,\dots,m,b_1,b_2,\dots,b_s)\eqno(3.11)$$
is true. The formulae $(3.5)$\,--\,$(3.7)$ and $(3.11)$ implies
that in the considered case the relations $(1.3)$ can be written
in the form $$d({\cal L}(f)^{n+1}_{(0,1,\dots,n)})={\cal
L}(f)^0_{(\,\,)}\p^{n+1}_{(0,1,\dots,n)}\,+$$ $$+\sum_{m=0}^{n-1}
\sum_{t=1}^{m+2}(-1)^{{\rm sign}(\sigma_{t,n-m})}{\cal
L}(f)^{m+1}_{(0,1,\dots,m)}\p^{n+1}_{(t-1,t,\dots,t+n-m-2)}-\p^{n+1}_{(0,1,\dots,n)}{\cal
L}(f)^{n+1}_{(\,\,)}\,-$$
$$-\sum_{m=0}^{n-1}\sum_{s=1}^{m+2}\sum_{N_{n-m}}\sum_{T_{m+2}}(-1)^{{\rm
sign}(\sigma_{t_1,\dots,t_s,n_1,\dots,n_s})}\p^{m+1}_{(0,1,\dots,m)}{\cal
L}(f)^{n+1}_{(b_1,b_2,\dots,b_s)},\eqno(3.12)$$ where by
$\sigma_{t,n-m}$ we denote the permutation $\sigma_{t_1,n_1}$ for
$t_1=t$ and $n_1=n-m$. Now, we compute ${\rm
sign}(\sigma_{t_1,\dots,t_s,n_1,\dots,n_s})$ for all $1\leqslant
s\leqslant m+2$. Denote by $|a_i|$ the number of elements in the
block $a_i$, where $1\leqslant i\leqslant s+1$, and by $|b_j|$ the
number of elements in the block $b_j$, where $1\leqslant
j\leqslant s$. Since $\sigma_{t_1,\dots,t_s,n_1,\dots,n_s})$ is a
permutation acting on the collection $(0,1,\dots,n)$ by
partitioning this collection into blocks
$(a_1,b_1,a_2,b_2,\dots,a_s,b_s,a_{s+1})$ and permuting of this
blocks by the rule $(3.8)$, the number of inversions
$I(\sigma_{t_1,\dots,t_s,n_1,\dots,n_s})$ of the permutation
$\sigma_{t_1,\dots,t_s,n_1,\dots,n_s})$ is equal
$$I(\sigma_{t_1,\dots,t_s,n_1,\dots,n_s})=$$
$$=|a_2||b_1|+|a_3|(|b_1|+|b_2|)+
\dots+|a_s|(|b_1|+\dots+|b_{s-1}|)+|a_{s+1}|(|b_1|+\dots+|b_s|).$$
By using the congruence
$I(\sigma_{t_1,\dots,t_s,n_1,\dots,n_s})\equiv {\rm
sign}(\sigma_{t_1,\dots,t_s,n_1,\dots,n_s})\,{\rm mod}(2)$ and the
equalities $$|a_i|=t_i,\quad 2\leqslant i\leqslant s,\quad
|a_{s+1}|=n+2-\sum_{k=1}^s(t_k+n_k),\quad \sum_{i=1}^s n_i=n-m,$$
we obtain the congruence $${\rm
sign}(\sigma_{t_1,\dots,t_s,n_1,\dots,n_s})\equiv
t_2n_1+t_3(n_1+n_2)+\dots+t_s(n_1+\dots+n_{s-1})\,+$$
$$+\,mn+m+(t_1+\dots+t_s)(n-m)\,{\rm mod}(2).$$ In particular, at
$s=1$, $t_1=t$, $n_1=n-m$ we have the congruence $${\rm
sign}(\sigma_{t,n-m})\equiv mn+m+t(n-m)\,{\rm mod}(2).$$ Now, we
show that the relations $(3.12)$ are equivalent to the relations
$(3.9)$. Indeed, by using $(3.5)$ and $(3.7)$ we write the
relations $(3.12)$ in the form $$d(f_{n+1})=f_0\pi_n
+\sum_{m=0}^{n-1}\sum_{t=1}^{m+2}(-1)^{\alpha}f_{m+1}(\underbrace{1\otimes\dots\otimes
1}_{t-1}\otimes\,\pi_{n-m-1}\otimes\underbrace{1\otimes\dots\otimes
1}_{m+2-t})\,-$$ $$-\,\pi_n(f_0\otimes\dots\otimes f_0)-
\sum_{m=0}^{n-1}\sum_{s=1}^{m+2}\sum_{N_{n-m}}\sum_{T_{m+2}}(-1)^\beta\pi_m(\underbrace{f_0\otimes\dots\otimes
f_0}_{t_1-1}\otimes \,f_{n_1}\otimes$$
$$\otimes\underbrace{f_0\otimes\dots\otimes
f_0}_{t_2-1}\otimes\,f_{n_2}\otimes\dots\otimes\underbrace{f_0\otimes\dots\otimes
f_0}_{t_s-1}\otimes
\,f_{n_s}\!\otimes\underbrace{f_0\otimes\dots\otimes
f_0}_{m+2-(t_1+\dots+t_s)}),$$ where $$\alpha=(n+1)(q-1)+{\rm
sign}(\sigma_{t,n-m})+(n-m)(q-1)\,+$$ $$+\,(m+1)(q+(n-m-1)-1),$$
$$\beta=(n+1)(q-1)+{\rm
sign}(\sigma_{t_1,\dots,t_s,n_1,\dots,n_s})+(n-m)(q-1)\,+$$
$$+\sum_{i=1}^{s-1}n_i(n_{i+1}+\dots+n_s)+(m+1)(q+(n-m)-1).$$ For
the exponent $\alpha$, we have $$\alpha\equiv
(n+1)(q-1)+mn+m+t(n-m)+(n-m)(q-1)\,+$$
$$+\,(m+1)(q+(n-m-1)-1)\equiv (m+1)(n-m-1)+mn+m+t(n-m)\equiv$$
$$\equiv t(n-m)+n+1\,{\rm mod}(2).$$ For the exponent $\beta$,
taking into account the equality $n-m=n_1+n_2+\dots+n_s$, we have
$$\beta\equiv
(n+1)(q-1)+t_2n_1+t_3(n_1+n_2)+\dots+t_s(n_1+\dots+n_{s-1})\,+$$
$$+\,mn+m+(t_1+\dots+t_s)(n-m)+(n-m)(q-1)+\sum_{i=1}^{s-1}n_i(n_{i+1}+\dots+n_s)\,+$$
$$+\,(m+1)(q+(n-m)-1)\equiv
t_2n_1+t_3(n_1+n_2)+\dots+t_s(n_1+\dots+n_{s-1})\,+$$
$$+\,n-m+(t_1+\dots+t_s)(n-m)+\sum_{i=1}^{s-1}n_i(n_{i+1}+\dots+n_s)\equiv$$
$$\equiv t_2n_1+t_3(n_1+n_2)+
\dots+t_s(n_1+\dots+n_{s-1})+(t_1-1)(n_1+\dots+n_s)\,+$$
$$+\,(t_2+\dots+t_s)(n_1+\dots+n_s)+
\sum_{i=1}^{s-1}n_i(n_{i+1}+\dots+n_s)\equiv$$ $$\equiv
\sum_{i=1}^s(t_i-1)(n_i+\dots+n_s)+\sum_{i=1}^{s-1}(n_i+1)(n_{i+1}+\dots+n_s)\,{\rm
mod}(2).$$ Thus, the relations $(3.12)$ are equivalent to the
relations $(3.9)$ and, consequently, the maps ${\cal
L}(f)^{n+1}_{(0,1,\dots,n)}=(-1)^{(n+1)(p-1)}f_{n+1}:(A^{\otimes
(n+2)})_p\to A'_{p+n+1}$, $n\geqslant 0$,  satisfy the relations
$(1.3)$. In a similar manner it is proved that the relations
$(1.3)$ holds for all defined by the formulas $(3.7)$ maps ${\cal
L}(f)^n_{(i_1,\dots,i_k)}$, where $i_k<n$.

Now, we check that the maps ${\cal L}(f)^n_{(i_1,\dots,i_k)}$,
where $i_k=n$, defined at $q=1$ and $s=0$ by $(3.8)$ satisfy the
relations $(1.3)$, i.e., we verify that the equality $$d({\cal
L}(f)^n_{(0,1,\dots,k-2,n)})=-\p^n_{(0,1,\dots,k-2,n)}{\cal
L}(f)^n_{(\,\,)}+{\cal
L}(f)^{n-k}_{(\,\,)}\p^n_{(0,1,\dots,k-2,n)}+$$
$$+\sum_{\varrho\in\Sigma_k}\sum_{I_{\varrho}}(-1)^{{\rm
sign}(\varrho)+1}\p^{n-k+m}_{(\widehat{\varrho(0)},\dots,\widehat{\varrho(m-1)})}
{\cal
L}(f)^n_{(\widehat{\varrho(m)},\dots,\widehat{\varrho(k-2)},\widehat{\varrho(n)})
}\,-$$ $$-\,{\cal
L}(f)^{n-k+m}_{(\widehat{\varrho(0)},\dots,\widehat{\varrho(m-1)})}
\p^n_{(\widehat{\varrho(m)},\dots,\widehat{\varrho(k-2)},\widehat{\varrho(n)})
}\eqno(3.13)$$ is true. By using the relations $dt_n=t_nd$, ${\cal
L}(f)^n_{(\,\,)}t_n=t_n{\cal L}(f)^n_{(\,\,)}$ and also the
conditions ${\cal L}(f)^n_{(0,1,\dots,k-2,n)}=(-1)^{k-1}{\cal
L}(f)^n_{(0,1,\dots,k-1)})t_n$,
$\p^n_{(0,1,\dots,k-2,n)}=(-1)^{k-1}\p^n_{(0,1,\dots,k-1)})t_n$,
we obtain $$d({\cal
L}(f)^n_{(0,1,\dots,k-2,n)})=-\p^n_{(0,1,\dots,k-2,n)}{\cal
L}(f)^n_{(\,\,)}+{\cal
L}(f)^{n-k}_{(\,\,)}\p^n_{(0,1,\dots,k-2,n)}\,+$$
$$+\sum_{\sigma\in\Sigma_k}\sum_{I_{\sigma}}(-1)^{{\rm
sign}(\sigma)+k}\p^{n-k+m}_{(\widehat{\sigma(0)},\dots,\widehat{\sigma(m-1)})}{\cal
L}(f)^n_{(\widehat{\sigma(m)},
\dots,\widehat{\sigma(k-1)})}t_n\,-$$ $$-\,{\cal
L}(f)^{n-k+m}_{(\widehat{\sigma(0)},\dots,\widehat{\sigma(m-1)})}
\p^n_{(\widehat{\sigma(m)},\dots,\widehat{\sigma(k-1)})}t_n.\eqno(3.14)$$
In the same way as it was done in the proof of Theorem 1.1, when
the coincidence of the right-hand sides of the equalities $(1.6)$
and $(1.7)$ was checked, it is proved that the right-hand sides of
the equalities $(3.13)$ and $(3.14)$ coincides. It follows that
the maps ${\cal L}(f)^n_{(0,1,\dots,k-2,n)}$ satisfy the relations
$(1.3)$. In similar way it is verified that the relations $(1.3)$
holds for all defined by $(3.8)$ maps ${\cal
L}(f)^n_{(i_1,\dots,i_k)}$, where $i_k=n$. Thus, all maps ${\cal
L}(f)^n_{(i_1,\dots,i_k)}\in{\cal L}(f)$ satisfy the the relations
$(1.3)$. Since above was shown that all maps ${\cal
L}(f)^n_{(i_1,\dots,i_k)}\in{\cal L}(f)$ satisfy the relations
$(1.5)$, the family of maps ${\cal L}(f)$ is a morphism of
$C\F$-modules ${\cal L}(f):({\cal L}(A),d,\p,t)\to ({\cal
L}(A'),d,\p,t)$.

Now, consider an arbitrary morphisms of $\A$-algebras
$f:(A,d,\pi_n)\to (A',d,\pi_n)$ and $g:(A',d,\pi_n)\to
(A'',d,\pi_n)$ and their composition $gf:(A,d,\pi_n)\to
(A'',d,\pi_n)$. We show that the equality of morphisms of
$C\F$-modules ${\cal L}(gf)={\cal L}(g){\cal L}(f)$ is true. We
must check that the maps ${\cal L}(gf)^n_{(i_1,\dots,i_k)}\in
{\cal L}(gf)$, $k\geqslant 0$, satisfy the relations $${\cal L}
(gf)^n_{(i_1,\dots,i_k)}=\sum_{\sigma\in\Sigma_k}\sum_{I'_{\sigma}}
(-1)^{{\rm sign}(\sigma)} {\cal
L}(g)^{n-k+m}_{(\widehat{\sigma(i_1)},\dots,\widehat{\sigma(i_m)})}
{\cal
L}(f)^n_{(\widehat{\sigma(i_{m+1})},\dots,\widehat{\sigma(i_k)})},\eqno(3.15)$$
where $I'_\sigma$ is the same as in $(1.4)$. Clearly, at $k=0$ we
have ${\cal L}(gf)^n_{(\,\,)}={\cal L}(g)^n_{(\,\,)}{\cal
L}(f)^n_{(\,\,)}$ because $(gf)_0^{\otimes (n+1)}=g_0^{\otimes
(n+1)}f_0^{\otimes (n+1)}$. Now, we check that the maps $${\cal
L}(gf)^{n+1}_{(0,1,\dots,n)}=(-1)^{(n+1)(p-1)}(gf)_{n+1}:(A^{\otimes
(n+2)})_p\to A''_{p+n+1},\quad n\geqslant 0,$$ satisfy the
relations $(3.15)$. With this purpose we write the relations
$(3.3)$ in the form
$$(gf)_{n+1}=g_0f_{n+1}+g_{n+1}(f_0\otimes\dots\otimes f_0)\,+$$
$$+\,\sum_{m=0}^{n-1}\sum_{s=1}^{m+2}\sum_{N_{n-m}}\sum_{T_{m+2}}(-1)^\mu
g_{m+1}(\underbrace{f_0\otimes\dots\otimes f_0}_{t_1-1}\otimes
\,f_{n_1}\otimes$$ $$\otimes\underbrace{f_0\otimes\dots\otimes
f_0}_{t_2-1}\otimes\,f_{n_2}\otimes\dots\otimes\underbrace{f_0\otimes\dots\otimes
f_0}_{t_s-1}\otimes
\,f_{n_s}\!\otimes\underbrace{f_0\otimes\dots\otimes
f_0}_{m+2-(t_1+\dots+t_s)}),\quad n\geqslant -1,\eqno(3.16)$$
where $N_{n-m}$ and $T_{m+2}$ and $\mu$ are the same as in
$(3.9)$. The formulae $(3.6)$ and $(3.7)$ follows that the
equalities $(3.16)$ can be written in the form $${\cal
L}(gf)^{n+1}_{(0,1,\dots,n)}={\cal L}(g)^0_{(\,\,)}{\cal
L}(f)^{n+1}_{(0,1,\dots,n)}+{\cal L}^{n+1}_{(0,1,\dots,n)}{\cal
L}(f)^{n+1}_{(\,\,)}+$$
$$+\sum_{m=0}^{n-1}\sum_{s=1}^{m+2}\sum_{N_{n-m}}\sum_{T_{m+2}}(-1)^\psi{\cal
L}^{m+1}_{(0,1,\dots,m)}{\cal L}(f)^{n+1}_{(b_1,b_2,\dots,b_s)},$$
where $\psi=\mu+(n+1)(q-1)+(m+1)(q+(n-m)-1)$ and number blocks
$b_1,\dots,b_s$ were specified above. On the other hand, the
formulae $(3.6)$, $(3.7)$ and $(3.11)$ follows that in the
considered case the relations $(3.15)$ can be written in the form
$${\cal L}(gf)^{n+1}_{(0,1,\dots,n)}={\cal L}(g)^0_{(\,\,)}{\cal
L}(f)^{n+1}_{(0,1,\dots,n)}+{\cal L}^{n+1}_{(0,1,\dots,n)}{\cal
L}(f)^{n+1}_{(\,\,)}+$$
$$+\sum_{m=0}^{n-1}\sum_{s=1}^{m+2}\sum_{N_{n-m}}\sum_{T_{m+2}}(-1)^{{\rm
sign}(\sigma_{t_1,\dots,t_s,n_1,\dots,n_s})}{\cal
L}^{m+1}_{(0,1,\dots,m)}{\cal L}(f)^{n+1}_{(b_1,b_2,\dots,b_s)},$$
where the permutation
$\sigma_{t_1,\dots,t_s,n_1,\dots,n_s}\in\Sigma_{n+1}$ is defined
by $(3.10)$. It was above shown that the congruence ${\rm
sign}(\sigma_{t_1,\dots,t_s,n_1,\dots,n_s})\equiv \psi\,{\rm
mod}(2)$ is true. Thus, the module maps ${\cal
L}(gf)^{n+1}_{(0,1,\dots,n)}=(-1)^{(n+1)(p-1)}(gf)_{n+1}:(A^{\otimes
(n+2)})_p\to A''_{p+n+1}$, $n\geqslant 0$, satisfy the relations
$(3.15)$. Similarly, it is proved that the relations $(3.14)$
holds for all maps ${\cal L}(gf)^n_{(i_1,\dots,i_k)}$, where
$i_k<n$. In the same way as it was done above, when we checked
that the maps ${\cal L}(f)^n_{(i_1,\dots,i_k)}$, where $i_k=n$,
satisfy the relations $(1.3)$, it is checked that the relations
$(3.15)$ holds for all maps ${\cal L}(gf)^n_{(i_1,\dots,i_k)}$,
where $i_k=n$. Thus, all maps ${\cal
L}(gf)^n_{(i_1,\dots,i_k)}\in{\cal L}(gf)$ satisfy the relations
$(3.15)$ and, consequently, the equality of morphisms of
$C\F$-modules ${\cal L}(gf)={\cal L}(g){\cal L}(f)$ is true.

The above considerations follows that there is the functor ${\cal
L}:\A(K)\to C\F(K)$. The required functor $HC:\A(K)\to GrM(K)$ we
define as a composition of the functor ${\cal L}:\A(K)\to C\F(K)$
and the functor $HC:C\F(K)\to GrM(K)$, which considered in Theorem
2.1.

Now, we show that the functor $HC:\A(K)\to GrM(K)$ sends homotopy
equivalences of $\A$-algebras into isomorphisms of graded modules.
Taking into account Theorem 2.1, it suffices to show that the
functor ${\cal L}:\A(K)\to C\F(K)$ sends homotopy equivalences
into homotopy equivalences of $C\F$-modules. With this purpose we
show that each homotopy between morphisms of $\A$-algebras induces
a homotopy between corresponding morphisms of $C\F$-modules.

Given any homotopy $h=\{h_n:(A^{\otimes(n+1)})_\bu\to
A'_{\bu+n+1}~|~n\in\mathbb{Z},~n\geqslant 0\}$ between morphisms
of $\A$-algebras $f:(A,d,\pi_n)\to (A',d,\pi_n)$ and
$g:(A,d,\pi_n)\to (A',d,\pi_n)$, we define a family of module maps
$${\cal L}(h)=\{{\cal L}(h)^n_{(i_1,\dots,i_k)}:{\cal
L}(A)_{n,p}\to {\cal L}(A')_{n-k,p+k+1}\},$$ $$n\geqslant 0,\quad
p\geqslant 0,\quad 0\leqslant k\leqslant n,\quad 0\leqslant
i_1<\dots<i_k\leqslant n,$$ by the following rules:

1$'$). If $k=0$, then $${\cal
L}(h)^n_{(\,\,)}=\sum\limits_{i=1}^{n+1}\underbrace{g_0\otimes\dots\otimes
g_0}_{i-1}\otimes\,h_0\otimes\underbrace{f_0\otimes\dots\otimes
f_0}_{n-i+1}\,;$$

2$'$). If $i_k<n$ and the collection
$$(i_1,\dots,i_k)=((j^1_1,\dots,j^1_{n_1}),(j^2_1,\dots,j^2_{n_2}),\dots,(j^s_1,\dots,j^s_{n_s}))$$
is the same as in the above rule 2) defining the formula $(3.7)$,
then $${\cal
L}(h)^n_{(i_1,\dots,i_k)}=(-1)^{k(p-1)+\gamma}\sum_{i=1}^s(-1)^{n_1+\dots+n_{i-1}}\underbrace{g_0\otimes\dots\otimes
g_0}_{k_1}\otimes\,g_{n_1}\!\otimes\dots$$
$$\dots\,\otimes\underbrace{g_0\otimes\dots\otimes
g_0}_{k_{i-1}}\otimes\,g_{n_{i-1}}\!\otimes\underbrace{g_0\otimes\dots\otimes
g_0}_{k_i}\otimes\,h_{n_i}\!\otimes\underbrace{f_0\otimes\dots\otimes
f_0}_{k_{i+1}}\otimes\,f_{n_{i+1}}\!\otimes\dots$$
$$\dots\otimes\underbrace{f_0\otimes\dots\otimes
f_0}_{k_s}\otimes\,f_{n_s}\!\otimes\underbrace{f_0\otimes\dots\otimes
f_0}_{k_{s+1}}\,+$$
$$+\,(-1)^{k(p-1)+\gamma}\sum_{i=1}^{s+1}(-1)^{n_1+\dots+n_{i-1}}\underbrace{g_0\otimes\dots\otimes
g_0}_{k_1}\otimes\,g_{n_1}\!\otimes\dots\otimes\underbrace{g_0\otimes\dots\otimes
g_0}_{k_{i-1}}\otimes\,g_{n_{i-1}}\otimes$$
$$\otimes\left\{\sum_{j=1}^{k_i}\underbrace{g_0\otimes\dots\otimes
g_0}_{j-1}\otimes\,h_0\otimes\underbrace{f_0\otimes\dots\otimes
f_0}_{k_i-j}\right\}\otimes\,f_{n_i}\!\otimes\underbrace{f_0\otimes\dots\otimes
f_0}_{k_{i+1}}\otimes\,f_{n_{i+1}}\otimes\dots$$
$$\dots\otimes\underbrace{f_0\otimes\dots\otimes
f_0}_{k_s}\otimes\,f_{n_s}\!\otimes\underbrace{f_0\otimes\dots\otimes
f_0}_{k_{s+1}}\,,$$ where $k_1,\dots,k_{s+1}$ and $\gamma$ are the
same as in $(3.7)$;

3$'$). If $i_k=n$ and the collection
$$(i_1,\dots,i_k)=((0,1,\dots,z-1-q),(j^1_1-q,\dots,j^1_{n_1}-q),(j^2_1-q,\dots,j^2_{n_2}-q),\dots$$
$$\dots,(j^s_1-q,\dots,j^s_{n_s}-q),(n-q+1,n-q+2,\dots,n))$$ is
the same as in the above rule 3) defining the formula $(3.8)$,
then
$${\cal L}(h)^n_{(i_1,\dots,i_k)}=(-1)^{q(z-1)}{\cal
L}(h)^n_{((0,1,\dots,z-1),(j^1_1,\dots,j^1_{n_1}),\dots,(j^s_1,
\dots,j^s_{n_s}))}t_n^q.$$

In similar way as it was done above in the case of morphisms of
$C\F$-modules ${\cal L}(f)$, it is proved that defined by any
homotopy $h=\{h_n:(A^{\otimes(n+1)})_\bu\to A'_{\bu+n+1}\}$
between morphisms of $\A$-algebras $f,g:(A,d,\pi_n)\to
(A',d,\pi_n)$ the family of maps ${\cal L}(h)=\{{\cal
L}(h)^n_{(i_1,\dots,i_k)}:{\cal L}(A)_{n,p}\to {\cal
L}(A')_{n-k,p+k+1}\}$ is a homotopy between morphisms of
$C\F$-modules ${\cal L}(f),{\cal L}(g):({\cal L}(A),d,\p,t)\to
({\cal L}(A'),d,\p,t)$. It follows that if $f:(A,d,\pi_n)\to
(A',d,\pi_n)$ is a homotopy equivalence of $\A$-algebras, then the
corresponding morphism ${\cal L}(f):({\cal L}(A),d,\p,t)\to ({\cal
L}(A'),d,\p,t)$ is a homotopy equivalence of $C\F$-modules. Thus,
the functor ${\cal L}:\A(K)\to C\F(K)$ sends homotopy equivalences
of $\A$-algebras into homotopy equivalences of $C\F$-modules and,
consequently, the functor $HC:\A(K)\to GrM(K)$ sends homotopy
equivalences of $\A$-algebras into isomorphisms of graded
modules.~~~$\blacksquare$

Let us consider applications of Theorem 3.1 to homology of
$\A$-algebras over any fields.

It is well known \cite{Kad} that if over any field given an
$\A$-algebra $(A,d,\pi_n)$, then on homologies $H(A)$ of this
$\A$-algebra, i.e., on homologies $H(A)$ of the chain complex
$(A,d)$, arises the $\A$-algebra structure $(H(A),d=0,\pi_n)$ such
that there is the homotopy equivalence of $\A$-algebras
$(A,d,\pi_n)\to (H(A),d=0,\pi_n)$. Applying Theorem 3.1 to this
situation, we obtain the following assertion.

{\bf Corollary 3.1}. The cyclic homology $HC(A)$ of any
$\A$-algebra $(A,d,\pi_n)$ over an arbitrary field is isomorphic
to the cyclic homology $HC(H(A))$ of the $\A$-algebra of
homologies $(H(A),d=0,\pi_n)$.~~~$\blacksquare$

In the case, when an $\A$-algebra $(A,d,\pi_n)$ is an associative
differential algebra $(A,d,\pi)$, where $\pi=\pi_0$ and $\pi_n=0$
at $n>0$, we have the following assertion.

{\bf Corollary 3.2}. The cyclic homology $HC(A)$ of any
associative differential algebra $(A,d,\pi)$ over an arbitrary
field is isomorphic to the cyclic homology $HC(H(A))$ of the
$\A$-algebra of homologies $(H(A),d=0,\pi_n)$.~~~$\blacksquare$

\vspace{1cm}

Serov Str., Saransk, Russia,

e-mail: slapin@mail.ru

\end{document}